\documentclass[12pt]{amsart}
\usepackage{latexsym}
\usepackage{amssymb, amsmath,mathrsfs,amsthm}
\usepackage{geometry}
\usepackage[OT2,T1]{fontenc}
\usepackage{setspace}
\usepackage{mathtools}
\usepackage{url}
\usepackage{enumerate}
\usepackage{float}
\usepackage{graphicx}

\counterwithin{figure}{section}

\newtheorem{theorem}{Theorem}[section]
\newtheorem*{theorem*}{Theorem}

\newtheorem{definition}[theorem]{Definition}

\newtheorem{corollary}[theorem]{Corollary}
\newtheorem{proposition}[theorem]{Proposition}

\usepackage{color}
\usepackage{forest}

\DeclareSymbolFont{cyrletters}{OT2}{wncyr}{m}{n}
\DeclareMathSymbol{\Sha}{\mathalpha}{cyrletters}{"58}

\theoremstyle{remark}
\newtheorem{remark}[theorem]{Remark}
\newtheorem{example}[theorem]{Example}

\begin{document}

 \title{Decomposable Blaschke products of degree $2^n$}

\begin{abstract}
We study the decomposability of a finite Blaschke product $B$ of degree $2^n$ into $n$ degree-$2$ Blaschke products, examining the connections between Blaschke products, the elliptical range theorem, Poncelet theorem, and the monodromy group. We show that if the numerical range of the compression of the shift operator, $W(S_B)$, with $B$ a Blaschke product of degree $n$, is an ellipse then $B$ can be written as a composition of lower-degree Blaschke products that correspond to a factorization of the integer $n$. We also show that a Blaschke product of degree $2^n$ with an elliptical Blaschke curve has at most $n$ distinct critical values, and we use this to examine the monodromy group associated with a regularized Blaschke product $B$. We prove that if $B$ can be decomposed into  $n$  degree-$2$ Blaschke products, then the monodromy group associated with $B$ is the wreath product of $n$ cyclic groups of order $2$. Lastly, we study the group of invariants of a Blaschke product $B$ of order $2^n$ when $B$ is a composition of $n$ Blaschke products of order $2$.

\end{abstract}
\author{Asuman G\"uven Aksoy}
\address{Department of Mathematical Sciences, Claremont McKenna College, Claremont, CA, 91711, USA.}
\email{asumanguvenaksoy@gmail.com}
\author{Francesca Arici}
\address{Mathematical Institute, Leiden University, P.O. Box 9512, 2300 RA Leiden, The Netherlands.}
\email{f.arici@math.leidenuniv.nl}
\author{M. Eugenia Celorrio}
\address{Department of Mathematics and Statistics, Flyde College, Lancaster University, Lancaster, LA1 4YF, United Kingdom.}
\email{celorrioramirez@lancaster.ac.uk}
\author{Pamela Gorkin}
\address{Department of Mathematics, Bucknell University, 377 Olin Science Building, Lewisburg, PA 17837, USA.}
\email{pgorkin@bucknell.edu}

\renewcommand{\thefootnote}{\fnsymbol{footnote}} 
\footnotetext{Primary 30J10; Secondary 30D05, 47A12, 20B05, 14N05}
\renewcommand{\thefootnote}{\arabic{footnote}}

\maketitle

\bigskip
\section{Introduction}

Blaschke products are important functions for the study of bounded analytic functions. They play the same role, hyperbolically speaking, on the unit disk $\mathbb{D}$ that polynomials play, in the Euclidean sense, on the plane, $\mathbb{C}$. Here we consider Blaschke products that are compositions of other, nontrivial, Blaschke products and consequences of this.

In the following, we give a brief review of the background, notation, and terminology that will be relevant to this paper. 
Let $\mathbb{D}$ denote the open unit disk and $\mathbb{T}$ denote the unit circle. A finite \textit{Blaschke product} $B$ of degree $n$ is a function of the form
$$B(z)= \gamma\prod_{j=1}^{n} \frac{z-a_j}{1-\overline{a_j}z},$$ 
where $a_j \in \mathbb{D}$ for $j=1,\dots, n$ and $\gamma\in \mathbb{T}$. Note that  Blaschke products of degree $1$ are  the disk automorphisms. Finite Blaschke products are $n$ to $1$  maps of the open unit disk $\mathbb{D}$ into itself and the unit circle $\mathbb{T}$  to itself. They are holomorphic on an open set containing the closed unit disk and have finitely many zeros in $\mathbb{D}$. We will consider the set of points in $\mathbb{T}$  that  the Blaschke product identifies; in other words, we will be interested in the solutions of $B(z)= \lambda$ for $\lambda \in \mathbb{T}$. Since the constant $\gamma$ will not play a role in the solution, we will take $\gamma=1$ in the above description of Blaschke product.  We say a Blaschke product $B$ is \textit{decomposable} if there exist Blaschke products $C$ and $D$ both of degree $n>1$ such that  $$ B(z)= C(D(z)) = (C\circ D)(z)$$ and $B$ is indecomposable otherwise. If the degree of $C$  and $D $ equals $k$ and $m$, respectively, then the degree of $B$ equals $km$. Observe that if $B$ is of prime degree, then $B$ is indecomposable.  Whether or not a Blachke product is a composition of two non-trivial Blaschke products has been studied  by many authors. The connection of decomposition to the monodromy groups are discussed in \cite{Ritt} for polynomials and \cite{Cowen} for Blaschke products.   In \cite{Wegert} and \cite{Beauty}, a visual representation of composition is discussed,  and in \cite{DGSSV15}  both algorithmic and geometric  arguments are presented, and the relationship between decomposable Blaschke products and curves with the Poncelet property are examined.  

Poncelet's theorem from projective geometry of conics has many deep consequences.  More than $200$ years ago  J. V. Poncelet discovered that if there exists a polygon of $n$-sides that is inscribed in a given conic and circumscribed about another conic, then infinitely many such polygons exist. This theorem is sometimes called Poncelet's porism and the related polygons are \textit{Poncelet's polygons}. The theorem, also called \textit{ Poncelet's closure theorem},  is a result about polygons that are inscribed in an ellipse and circumscribe a smaller ellipse and has been studied in several settings (see, for example, \cite{DelCentina, Dragovic2,  Flatto, HalbeisenHungerbuhler, Mirman}, among others). Later,  Darboux  found a new proof of the Poncelet closure theorem based on the properties of certain curves, known as  \textit{Poncelet-Darboux curves}, \cite{Dragovic}.  These are the curves of degree $n$  passing through the intersection points of $n+1$  tangents to a given conic. When we say  $C$ is an \textit{$n$-Poncelet curve} in $\mathbb{D}$,  we mean for every $\lambda$ with $|\lambda|=1$, there exists an $n$-sided convex polygon that circumscribes $C$ that has $\lambda$ as one vertex and all vertices on $\mathbb{T}$.

Poncelet curves and Blaschke products are also related to a certain class of square matrices via the numerical range of a matrix. Recall that given an $n \times n  $ matrix $A$, the \textit{numerical range}  $W(A)$ is defined by
$$W(A):= \left\{ \langle Ax,x \rangle: \,\,x\in \mathbb{C}^n,\,\, ||x|| =1\right\}.$$ The set $W(A)$ contains the spectrum of $A$, it is a convex set, and its outer boundary is a convex curve.  In particular if $A$ is a $2 \times 2$ matrix, then its numerical range is either a point, a line segment, or an elliptical disk -- all of these can be thought of as elliptical disks.  This is called the \textit{elliptical range theorem} \cite{ DGSV18}. Although this is a theorem about  $2\times 2$ matrices, it also sheds light on the numerical range of $n \times n$ matrices. Surprisingly, Blaschke products, the numerical range and Poncelet's theorem are connected. Roughly speaking, this connection can be described by noting that the convex hull of each of the circumscribing polygons with vertices on $\mathbb{T}$ represents the numerical range of a certain unitary matrix that is related to a Blaschke product via an operator that is a {\it compression of the shift operator}. See, for example,  \cite{DGSV18} for an overview of these connections and Section~\ref{sec:Ellipses, Numerical Range, and the Blaschke curve} below.  

In  \cite{Fujimura, Fujimura2}, Fujimura considered geometric properties of Blaschke products of degree $n$ and the line segments that are tangent to the unit circle at the points $B$ maps to $\lambda$ on the unit circle as well as the line segments joining successive points. For degree-$4$ Blaschke products, she considered Blaschke products for which the trace of these lines is an ellipse. Using Gr\"obner bases, she showed that an ellipse is inscribed in a quadrilateral that is inscribed in $\mathbb{T}$ if and only if the Blaschke product is a composition of two degree-$2$ Blaschke products. In \cite{GorkinWagner}, the authors gave an operator-theoretic proof of this result. They showed that an ellipse is a Poncelet ellipse (in the sense described above) inscribed in a quadrilateral inscribed in $\mathbb{T}$ if and only if the ellipse is the Blaschke curve associated with a Blaschke product $\widehat{B}(z) = z B(z)$ and the compression of the shift operator $S_B$ (see Section~\ref{sec:compression}) associated with the Blaschke product $B$ has elliptical numerical range. Further, $S_B$ has elliptical numerical range if and only if $\widehat{B}$ is a composition of two degree-$2$ Blaschke products. In  \cite{Simanek3} the authors prove obtain a similar result for Blaschke products of degree $6$.  

In this paper, we consider a Blaschke product $B$ of degree $n=2^k -1$ and $\widehat{B}$. For each $\lambda \in \mathbb{T}$ we let $z_{1, \lambda}, \ldots, z_{n+1, \lambda}$ denote the points on $\mathbb{T}$ that $\hat{B}$ maps to $B(\lambda)$ ordered by increasing argument. By \cite{GauWu, Mirman} each convex polygon connecting these points circumscribes a convex smooth curve that we call the {\it Blaschke curve}. In Section~\ref{sec:Ellipses, Numerical Range, and the Blaschke curve} we show that if the numerical range, $W(S_B)$, is an elliptical disk, then every lower degree curve in  the Poncelet package (see that section for a definition) for $\widehat{B}(z) := zB(z)$ is an ellipse and $\widehat{B}$ is a composition of $k$ degree-$2$ Blaschke products.  Here, the boundary of the numerical range is the Blaschke curve associated with $\widehat{B}$. This statement about composition can be generalized to Blaschke products of degree $n$; see Theorem~\ref{thm:generaln} below. We also give an example of  a Blaschke product of degree-$8$ for which  $W(S_B)$ is  elliptical and a ``non-example'' of a Blaschke product of degree-$8$ that factors into three degree-$2$ Blaschke products, but such that $W(S_B)$ is not  an elliptical disk. 
 
In Section~\ref{sec:critical values}, we turn our attention to a deep theorem of Ritt \cite {Ritt}  that classifies decomposability  of $B$ in terms of the monodromy group associated with $B$.  See also \cite{Cowen, Wegert}  for more recent developments. For this purpose, we start by examining critical values of Blaschke products with elliptical Blaschke curve. As usual, by the set of critical values of $B$ we mean  $\{ w\in \mathbb{D}: \,\, w=B(z) \,\,\mbox{and}\,\, B'(z)=0\}.$ 
We prove that if one writes a Blaschke product $B$ of degree $2^n$ as composition of  $n$ degree-$2$ Blaschke products, then $B$ has at most $n$ distinct critical values. After obtaining a description of a Blaschke product of degree-$n$ with one critical value,  we deduce that if $B$ is such a Blaschke product of degree $n= p_1p_2\dots p_m$, then $B$ can be factored in any order as a composition of $m$ Blaschke products of  degree $n= p_1,p_2, \dots ,p_m$ .
  
In Section~\ref{sec:monodromy}, we describe the monodromy group of $B$ assuming that $B$ is a normalized Blaschke product (see Section~\ref{sec:monodromy} for this definition) that is the composition of $n$ degree-$2$ Blaschke products with $n$ critical values.  We prove that in this case the monodromy group associated with $B$ is the wreath product of $n$ cyclic groups of order $2$. 

In Section~\ref{sec:group of invariants}, letting $C(\mathbb{T})$ denote the space of continuous functions from the unit circle to itself, we study the {\it group of invariants} of a finite Blaschke products $B$; that is, the group \[\mathcal{G}_B = \{u \in C(\mathbb{T}): B \circ u = B\}.\] In \cite{CassierChalendar} Cassier and Chalendar showed that the group of invariants of a Blaschke product of degree $n$ is a cyclic group of order $n$. The group of invariants for infinite Blaschke products with finitely many singularities was considered in \cite{CGP12}. Here we show that if $B$ is a composition of $n$ degree-$2$ Blaschke products, say $B = C_n \circ C_{n-1} \circ \ldots \circ C_1$, then the group of invariants of $C_{j} \circ \cdots \circ C_1$ is of index $2$ (and hence normal) in the group of invariants of $C_{j + 1} \circ \cdots \circ C_1$. From this, we are able to obtain a connection between elements of the group and a particular automorphism of the unit disk of the form $\varphi_a(z) = (z-a)/(1 - \overline{a} z)$ for $a \in \mathbb{D}$.

\section{Closure results}
\raggedbottom
 We begin by discussing the background required in real projective space $\mathbb{P}^2(\mathbb{K})$, where $\mathbb{K}$ is a field. In this paper, our field will be $\mathbb{R}$ or $\mathbb{C}$. 
 For $x, y, z \in \mathbb{K}$, we have an equivalence relation defined on $\mathbb{K}^3 \setminus \{(0,0,0)\}$ by $(x, y, z) = (x^\prime, y^\prime, z^\prime)$ if there exists a scalar $\lambda \ne 0$ such that $x^\prime = \lambda x, y^\prime = \lambda y,$ and $z^\prime = \lambda z$.  Points in the real (complex) projective plane $\mathbb{P}^2(\mathbb{R})$ are equivalence classes of triples of real (resp. complex) numbers for the relation above. The real projective plane $\mathbb{P}^2(\mathbb{R})$ is embedded in the complex projective plane $\mathbb{P}^2(\mathbb{C})$. 

An algebraic curve in $\mathbb{P}^2(\mathbb{C})$ is the set of zeros of a homogeneous polynomial $f$ with complex coefficients. 
A \emph{real} algebraic curve is an algebraic curve in the complex projective plane given by an equation $f(x, y, z) = 0$, where $f$ is a homogeneous polynomial with real coefficients.  The set of real points of a curve $\Gamma$ will be denoted by $\Gamma_\mathbb{R}$ and is defined by $\Gamma_\mathbb{R} = \Gamma \cap \mathbb{P}^2(\mathbb{R})$.

Since polynomial rings over a field are unique factorization domains, every algebraic curve $C$ is the union of finitely many irreducible curves, called its irreducible components.  If $C_1, \ldots, C_k$ are the irreducible components of  $C$ with irreducible defining polynomials $f_1,\ldots,f_k$, then $f = f_1 \cdots f_k$ is the minimal polynomial defining $C$. 

\subsection{Duality and reciprocation about $\mathbb{T}$}

If $\Gamma$ is a real algebraic curve in $\mathbb{P}^2(\mathbb{C})$, then the dual of $\Gamma$ is denoted by $\Gamma^\ast$, where the points correspond to the tangent lines to $\Gamma$.

Note that the dual curve of a general plane algebraic curve $\Gamma$ is the union of the dual curves of its irreducible components. In particular, $\Gamma$ and $\Gamma^*$ have the same number of irreducible components.
Moreover, note that the dual to a conic is a conic and the dual to a line is a point. By \cite[Theorem~5.1]{Fis}, duality is involutive, that is $(\Gamma^*)^* = \Gamma$. Note that the result is stated under the assumption that the curve has no lines as components, because a point is not an algebraic curve. 


\subsection*{Poncelet's porism}
Having now established some notation, we turn to some relevant theorems and to various formulations thereof. Poncelet discovered that if there exists a polygon of $n$-sides that is inscribed in a given conic and circumscribed about another conic, then infinitely many such polygons exist. This theorem is sometimes called Poncelet's porism and the related polygons were called Poncelet's polygons. The two conics are said to be $n$-Poncelet related. 

Suppose that an ellipse $E_1$ is inscribed in a convex $n$-gon that is itself inscribed in the unit circle; that is, $E_1$ and the unit circle are Poncelet related.  Consider the diagonals that connect vertex $k - 1$ to vertex $k + m$, for $k = 1, \ldots, n$.  Then \cite [p. 208]{Berger} the envelope of these diagonals is again an ellipse, $E_{m + 1}$, for $m = 1, \ldots, [n/2]-1$, where $[x]$ denotes the greatest integer less than or equal to $x$. Below are different versions of the theorems that we use later. These are often referred to as Darboux's theorem. The one we will use can be found in \cite{Mirman2} or \cite[p. 208]{Berger}.


\begin{theorem} \label{thm:darboux2}
\label{Darboux} Let $E_1$ be an ellipse inscribed in a convex $n$-gon that is, in turn, inscribed in $\mathbb{T}$. Consider the diagonals of the Poncelet polygons that leap over $m$ vertices (i.e., all such diagonals that connect vertex $k - 1$ with vertex $k + m$ where $k = 1, 2, \ldots, n$ and the vertex numbers are taken modulo $n$). The envelope of these diagonals is also an ellipse $E_{m + 1}$ for $m = 1, \ldots, [n/2]-1$.\end{theorem}

These curves were also studied in \cite{Dragovic}, where the author looks at Poncelet--Darboux theorems. There, the relevant result relies on the notion of a conic and a curve being Poncelet--Darboux related, an extension of the notion for two conics being Poncelet related.

The result from Darboux's theorem,  and the notion of being Poncelet related, can be extended to the setting of a conic $C$ and a curve $S$ of degree $n-1$ (cf. \cite[Subsection 2.1]{Dragovic}).

\begin{definition}
Let $C$ be a conic and  $S$ a curve of degree $n-1$.  We say that $S$ and $C$ are Poncelet--Darboux $n$-related,  if $S$ contains all the intersection points of the n-tangents to $C$.  
\end{definition}
Note that the above definition implicitly includes the main result of a theorem by Darboux (\cite{Darboux}, \cite[Theorem 1]{Dragovic}). 
 If $S$ contains all the intersection points of $n$ given tangents,  then it will contain the intersection points of any other set of $n$-tangents. 

\begin{theorem}[Darboux's Theorem,  cf. {\cite[Theorem 5]{Dragovic}}] 
\label{thm:Dragovic} If a curve $S$ of degree $n-1$ is $n$-Poncelet-Darboux related to a conic $K$ and if there is a conic $C$, a component of $S$ which is $n$-Poncelet related to the conic $K$, then for $n$ odd, the curve $S$ can be completely decomposed into $(n-1)/2$
 conics, and for $n$ even, it can be decomposed into $(n-2)/2$ conics and a line.
  \end{theorem}

In other words,  the minimal homogeneous polynomial defining the curve $S$ factors as a product of 
$(n-1)/2$ degree-2 irreducible polynomials if $n$ is odd,  and $(n-2)/2$ degree-2 irreducible polynomials and a degree one polynomial if $n$ is even. 

In the case where $K$ is an ellipse,  and hence also for $K=\mathbb{T}$,  the tangency condition in the theorem above implies that $S$ lies outside $K$.  As we are interested in studying curves inside the unit circle,  we need a dual formulation to Theorem~\ref{thm:Dragovic}. This can be obtained by considering $C=\mathbb{T}$ and $S$ as the ``dual curve'' to a closed convex curve in $\mathbb{D}$. By reciprocation, $S$ will lie on the outside $\mathbb{T}$.

\begin{theorem}[{cf. \cite[Theorem B]{Simanek2}}]\label{thm:Simanek}
Let $C$ be a closed convex curve in $\mathbb{D}$ and suppose that
there is an $n$-sided polygon inscribed in $\mathbb{T}$ and circumscribed about $C$.  Assume further that $C$
is a connected component of a real algebraic curve $\Gamma$ in $\mathbb{D}$ of class $n- 1$ such that each
diagonal of the polygon is tangent to $\Gamma$.  Then for every point $\lambda$ of $\mathbb{T}$ there is an $n$-sided convex polygon that is inscribed in $\mathbb{T}$, circumscribed about $\Gamma$, and has $\lambda$ as a vertex. 
In the special case when $C$ is an ellipse,  the curve $\Gamma$ decomposes into $(n - 1)/2$ ellipses if $n$ is
odd, and $(n -2)/2$ ellipses and an isolated point if $n$ is even.
\end{theorem}


%


Motivated by the above formulation, and in accordance with Mirman \cite{Mirman2}, we will say that a smooth closed curve $\Gamma$ contained in $\mathbb{D}$ is an $n$-Poncelet curve if for every point $\lambda$ of $\mathbb{T}$ there is an $n$-sided convex polygon that is inscribed in $\mathbb{T}$, circumscribed about $\Gamma$, and has $\lambda$ as a vertex. 

 If we begin with a Poncelet curve that is inscribed in a convex $n$-gon, Mirman considers the diagonals of these polygons and denotes the envelope of the diagonals that skip $m$ vertices, with $m \le [n/2]-1$, by $K_{m + 1}$. The set $K_1, \ldots, K_{[n/2]}$ is called a package of Poncelet curves. Two recent papers \cite{Simanek2}, \cite{Simanek3} provide many details and examples relevant to Mirman's work as well as this paper. In our setting, the $K_j$ will be ellipses and for $j > 1$ we sometimes refer to these as lower-degree curves. Thus, 
$\bigcup_{j = 1}^{[n/2]} K_j$ is the package of Poncelet curves generated by the convex $n$-gons. 

\subsection{Poncelet, Darboux and the numerical range}

We will be interested in this in the following setting: Given an $n \times n$ matrix $A$ we let the real part of $A$ and the imaginary part be the self-adjoint matrices defined by 
\[\mbox{Re}(A) = \frac{A + A^\ast}{2} \qquad \mbox{and} \qquad \mbox{Im}(A) = \frac{A - A^\ast}{2 i}.\] Of course, $A = \mbox{Re}(A) + i \mbox{Im}(A)$. By \cite{Kippenhahn}, we may associate a curve $\Gamma$ of class $n$ in homogeneous line coordinates via the function $f_A(u, v, w) = ~\mbox{det}(u \, \mbox{Re}(A) + v  \, \mbox{Im}(A) + w \, I)$.  Consider the algebraic curve $C(A)$ determined by $f_A = 0$ in $\mathbb{P}^2(\mathbb{C})$; that is
\[C(A) = \{(u, v, w) \in \mathbb{P}^2(\mathbb{C}): u \, x + v \, y + w \, z = 0~\mbox{is a tangent line to}~ f_A(x, y, z) = 0\}.\] Kippenhahn's theorem says that $W(A)$ is the convex hull of the real points of $C(A)$.  This curve, $C_{\mathbb{R}}(A)$, is called the Kippenhahn curve of $A$.

We will be interested in the numerical range of compressions of the shift operator associated with finite Blaschke products.  If we consider $S_B$, the vertices of each polygon are determined by the function $\widehat{B}(z) := z B(z)$ as follows: Given $\lambda \in \mathbb{T}$, the vertices of the corresponding polygon are the solutions of $\widehat{B}(z) = \widehat{B}(\lambda)$. Such operators have no unitary summand, so by \cite[Theorem 2]{Kippenhahn} their eigenvalues are interior to $W(A)$ and the boundary of $W(A)$ is smooth (see, for example, \cite[Theorem 12]{Kippenhahn}). Our curves have the property that for $\lambda \in \mathbb{T}$ the two sides of the polygons with vertex at $\lambda$ are tangent to $\Gamma$, and every point of $\Gamma$ is such a point of tangency, \cite{GauWu}. We see from the expression for $f_A$ that $f_A$ is a homogeneous polynomial of degree $n$ with real coefficients, which tells us that the dual of $\Gamma$, denoted by $\Gamma^\ast$, is a real algebraic curve in $\mathbb{P}^2(\mathbb{C})$.  It is known \cite[Lemma 3.10]{Simanek2} that if $\Gamma$ is a real algebraic curve of class $n-1$, then 
\[\bigcup_{j = 1}^{[n/2]} K_j = \Gamma.\]





\section{Ellipses, Numerical Range, and the Blaschke curve}\label{sec:Ellipses, Numerical Range, and the Blaschke curve}

In \cite{DaeppGorkinMortini}  it is shown that for an ellipse $E$, the following conditions are equivalent (see also \cite{Fujimura}):
\begin{enumerate}
\item There exists a triangle that $E$ is inscribed in and that the unit circle is circumscribed about.
\item For some $a,b$  in the unit disk, $E$ is defined by the equation
$$ |z-a|+|z-b| = |1- \overline{a}b|.$$
\end{enumerate}
Fujimura extended this result to degree-$4$  Blaschke products, giving  necessary and sufficient conditions for a quadrilateral inscribed in $\mathbb{T}$ to circumscribe an ellipse. Specifically, an ellipse $E$ is inscribed in a quadrilateral  that is itself inscribed in the unit circle if and only if $E$ is associated with a Blaschke product  $B = C\circ D$,  where $C$ and $D$ are  degree-$2$ Blaschke products. Later Gorkin and Wagner in  \cite{GorkinWagner} put this result into its operator theoretic context by characterizing  the numerical range of the compression of the shift operator $S_B$ and giving  the placement of the three zeros of $B$. In fact they show that the location of two of the three zeros of $B$  completely determines the location of the third. They prove that for any  two points $b$ and $c$ in the unit disk there is a uniquie Poncelet $4$-ellipse with foci at $b$ and $c$, and that the corresponding Blaschke product has zeros at $0,b,c$ and $a$, where $a$ is given in terms of $b$ and $c$. Given a Blaschke product $B$ of degree $n$ and a point $\lambda \in \mathbb{T}$, 

In \cite{Fujimura2}, Fujimura considered the set of lines joining two distinct preimages in $B^{-1}(\lambda)$ and the envelope of these lines, which she called the interior curve associated with the Blaschke product. In addition, if one considers the set of lines tangent to the unit circle, the trace of the intersection of each of two lines in this set is called the exterior curve. Several geometric properties about these two curves are presented. In particular, \cite[Theorem 5]{Fujimura2} contains a relation between the dual of the homogenized exterior curve and the interior curve.

\subsection{Blaschke products and composition}
\label{sec:compression}

We begin by considering the Blaschke curve, which is defined as follows.

 \begin{definition}
Let $B$ be a Blaschke product of degree $n+1$. Then the Blaschke curve  $\mathcal {C}$ associated with this  Blaschke product  is a  curve inscribed in the convex polygons  with vertices at the solutions $z_j \in \mathbb{T}$ of  $B(z_j)=\lambda$  for each $\lambda \in \mathbb{T}$. Each point on the curve is the point of tangency of such a circumscribing convex polygon. \end{definition}

Note that a Blaschke curve is assumed to be contained in $\mathbb{D}$. It is known that the line joining the  $z_j$ to $z_{j+1}$  is tangent to $\mathcal{C}$ at a single point, the curve contains no line segments, is a differentiable algebraic curve, and every point on the curve $\mathcal{C}$ can be obtained using the Blaschke product as described in the definition. See \cite{GauWu, GauWu2, DaeppGorkinVoss}.

Note that if $\varphi$ is an automorphism of the disk and $B$ is a Blaschke  product, then $\varphi \circ B$ and $B$ have the same Blaschke curve, but in general different Blaschke products produce different  Blaschke curves. Furthermore, even though a Blaschke curve is a Poncelet curve, not all Poncelet curves are Blaschke curves, \cite{Mirman98}. The relationship between the geometry  of the numerical range (elliptical numerical range) and  composition of degree two Blaschke products given in [10]. When we assume that $B(0) = 0$, writing $B(z) = z B_1(z)$, it can be shown that the Blaschke curve is the smooth convex curve in $\mathbb{D}$ that is the boundary of the an operator that is unitarily equivalent to the compressed shift operator $S_{B_1} = P_{B_1} S|_{K_{B_1}}$, where $H^2$ is the Hardy space, $K_{B_1} = H^2 \ominus B_1 H^2$ is the model space, $S$ is the shift operator, and $P_{B_1}$ is the orthogonal projection of $H^2$ onto $K_{B_1}$ (see \cite{DGSV18, GauWu, GauWu2}.) This class consists of contractions $T$ with eigenvalues inside the open unit disk $\mathbb{D}$, that satisfy $\mbox{rank}~(I - T^\star T) = 1$.
 
 The proof that an elliptical numerical range implies that Blaschke products are compositions is an extension of theorems for degree-$4$ Blaschke products \cite{Fujimura, GorkinWagner}, and degree-$6$ Blaschke products in \cite{Simanek3}. We will use the following theorem:
 
 \begin{theorem}\cite[Theorem 2.3]{DGSSV15}\label{thm:unique} Given two sets of points $z_1, \ldots, z_n$ and $z_1^\prime, \ldots, z_n^\prime$ interlaced on the unit circle, there is a Blaschke product $B$ of degree $n$ such that $B(0) = 0$, $B(z_j) = B(z_k)$ and $B(z_j^\prime) = B(z_k^\prime)$ for all $j$ and $k$. This Blaschke product $B$ is unique up to a rotation factor $\lambda$ with $|\lambda| = 1$. \end{theorem}
 
 \begin{theorem}\label{thm:2} Let $B$ be a Blaschke product of degree $2^n-1$. If $W(S_B)$ is an elliptical disk then every lower degree curve in the Poncelet package for $\widehat{B}(z):=z B(z)$ is an ellipse or a point and $\widehat{B}$ is a composition of $k$ degree-$2$ Blaschke products. \end{theorem}
 
\begin{proof} By \cite[Theorem 2.1]{GauWu}, the boundary of $W(S_B)$ is circumscribed by a $2^n$-sided convex polygon. \color{black} Applying Darboux's theorem as given in Theorem~\ref{Darboux} implies, among other things, that all curves inscribed in the appropriate convex polygons with $2^m$ sides, $m = 2, 3, \ldots, n$ are elliptical or a point. In addition, by Theorem~\ref{thm:Simanek}, since the degree of the dual curve is odd, the algebraic curve in question, which is the dual of the dual, will decompose into conics and a point. The set of diagonals  joining vertex $k$ with vertex $k + 2^{n-1}$ (mod $2^n$) will yield the point. It is shown in \cite{GorkinWagner} that if $\widehat{B}$ has degree $4$ and elliptical numerical range then $\widehat{B}$ is the composition of two degree-$2$ Blaschke products. We prove the rest by induction. 

So suppose that if $\widehat{B_1}$, defined by $\widehat{B_1}(z) = z B_1(z)$,  is a Blaschke product with corresponding Blaschke curve (that is, the boundary of $W(S_{B_1})$) elliptical and the degree of $\widehat{B}_1$ is equal to $2^{n_1}$ with $2 \le n_1 < n$, then $\widehat{B_1}$ is the composition of $n_1$ degree-$2$ Blaschke products. Now consider $\widehat{B}$ of degree $2^n$. According to Theorem~\ref{Darboux}, if we have an ellipse that is inscribed in a convex $n$-gon that is itself inscribed in $\mathbb{T}$, then the diagonals of the circumscribing polygons (that connect vertex $k-1$ with vertex $k + m$ for $k = 1, 2, \ldots, 2^n$ and indices chosen modulo $2^n$), have as envelope an ellipse $E_m$.  Take the vertices of two polygons, $P_z$ and $P_w$, and denote them by $\{z_0, \ldots, z_{2^n-1}\}$ and $\{w_0, \ldots, w_{2^n-1}\}$. If we skip a point when connecting points identified by $\widehat{B}$ on $\mathbb{T}$, Theorem~\ref{Darboux} tells us that we will see a Poncelet ellipse and since these polygons have $2^{n-1}$ interlaced vertices, these line segments will produce a closed convex polygon. So, skipping a point in each set $\{z_j\}$ and $\{w_j\}$, we obtain four convex polygons with $2^{n-1}$ vertices $\{z_{2j}\}$, $\{z_{2j+1}\}$, $\{w_{2j}\}$, and $\{w_{2j+1}\}$. By Theorem~\ref{Darboux}, there is one ellipse inscribed in these convex polygons. That ellipse is a Poncelet ellipse and therefore, by \cite[p. 219]{GauWu3} it is also a Blaschke curve. Thus, there is a Blaschke product $D$ of degree $2^{n-1}$ with $D(0) = 0$ that identifies each point in a set with every other point in the same set: For $j \ne l$,
 \[D(z_{2j}) = D(z_{2l}) = \lambda_1, D(z_{2j+1}) = D(z_{2l+1}) = \lambda_2, \] and
 \[D(w_{2j}) = D(w_{2l}) = \gamma_1, ~\mbox{and}~ D(w_{2j+1}) = D(w_{2l + 1}) = \gamma_2.\] Now because the points are interlaced and a Blaschke product has increasing argument on the unit circle, we may assume without loss of generality that
 \[\mbox{arg}(\lambda_1) < \mbox{arg}(\gamma_1) < \mbox{arg}(\lambda_2) < \mbox{arg}(\gamma_2).\]
 Therefore, by \cite[Theorem 9]{GorkinRhoades} there is a degree-$2$ Blaschke product $C$ mapping $0$ to $0$ such that
 \[C(\lambda_1) = C(\lambda_2)~\mbox{and}~C(\gamma_1) = C(\gamma_2).\] Thus $C \circ D$ identifies the vertices of $P_z$ and $C \circ D$ identifies the vertices of $P_w$. Further $C \circ D(0) = 0$. 
 


By the uniqueness guaranteed by Theorem~\ref{thm:unique}, there exists $\lambda \in \mathbb{T}$ such that $\widehat{B} = \lambda (C \circ D)$. Since $D$ is degree $2^{n-1}$ with $D(0) = 0$ that identifies every other point in our sets, Theorem~\ref{Darboux} applies to $D$ and therefore the induction hypothesis applies to $D$. Thus $D$ factors into a composition of $n-1$ Blaschke products of degree $2$ and therefore the result holds. \end{proof}

\begin{remark} If all the lower-degree curves (that is, when we skip at least one vertex) are ellipses, one can use a counting argument to show that the curve that we obtain by connecting successive points is also a conic: Consider the boundary of the numerical range. By Kippenhahn's theorem \cite[Theorem 10]{Kippenhahn} this is an algebraic curve and the dual curve $\Gamma$ has degree $2^n - 1$. Since the degree of the dual curve is odd, and we assume that each lower-degree component is an ellipse, no component will be a point. Given our assumptions, the numerical range is contained in the open unit disk and the lower-degree curves are all ellipses. Therefore, since the ellipses skipping more than $2^{n-1}$ points can be matched with one of those skipping fewer than $2^{n-1}$ points and the one skipping exactly $2^{n-1}$ points yields a single point, we get $2^{n-1} - 2$ ellipses for the lower-degree cases. (For example, in case we have $8 = 2^3$, we have ellipses when we join every $2$nd or $3$rd point, and we get a point when we join every $4$th point. The case when we join every $5$th point is the same as joining every $3$rd point and joining every $6$th point is the same as joining every second point. The case when we join subsequent points is the one we are trying to determine.) There are $2^{n-1}-2$ ellipses, a line, and the curve we are trying to identify. Now the dual of an algebraic curve of degree $2$ maintains the same degree and the dual of a line is a point, so the degrees of the dual curves corresponding to components that we obtain by skipping at least one point therefore total $2(2^{n-1} - 2) = 2^n - 4$ for those corresponding to ellipses and $1$ for the line, or $2^n -3$. But we should have degree $2^n-1$, so the component of $\Gamma$ that we have not yet counted, namely the one corresponding to the curve in which we do not skip any points, must have degree $2$. Therefore, it must be a conic and the dual of the dual (the original curve) is contained in the boundary of the numerical range and a compact convex subset of $\mathbb{D}$. Since the dual of a degree-$2$ curve is a conic, it must be an ellipse.  In Example ~\ref{ex:notelliptical} we present a Blaschke product of degree $2^n$ such that $W(S_B)$ is not elliptical, but all lower-degree curves corresponding to polygons that are Poncelet curves inscribed in convex polygons with $2^m$ sides with $m < n$ are elliptical. Of course, Theorem~\ref{thm:Dragovic} tells us that none of the curves that are $2^n$-Poncelet can be elliptical. This is also illustrated in this example. \end{remark}


The proof of Theorem~\ref{thm:2} works in more generality, as indicated below. We have stated it in this way because of our focus on Blaschke products that have degree a power of $2$.

\begin{theorem}\label{thm:generaln} Let $B$ be a degree $n$ Blaschke product with an elliptical Blaschke curve. Then for each factor $k > 1$ of $n$, there are  Blaschke products $C$ of degree $k$ and $D$ of degree $n/k$ such that $B = C \circ D$.
\end{theorem}

The proof is essentially the same as that of Theorem~\ref{thm:2} above. We provide a brief outline of the proof, indicating places where the proof will be slightly different.

\begin{proof} Suppose $n = k m$, with $k, m \in \mathbb{N} \setminus \{1\}$. Let $P_1$ and $P_2$ denote two Poncelet polygons with $n$ vertices, $z_1, \ldots, z_n$ and $w_1, \ldots, w_n$. Using every $k$-th point as a vertex, we get $k$ convex $m$-gons and, applying Theorem~\ref{Darboux}, we may conclude that  they circumscribe the same ellipse. Since this ellipse is a Poncelet curve contained in $\mathbb{D}$ inscribed in a convex polygon that has all of its vertices on $\mathbb{T}$, as above there is a Blaschke product $D$ of degree $m$ that maps $0$ to $0$ and identifies each set of $m$ vertices of each of the respective polygons. That is, these $m$-gons are the convex polygons circumscribing the Blaschke curve of $D$. Now there are $k$ polygons with $m$ vertices and $D$ is exactly $m$-to-$1$, so $D$ takes $k$ values on the $k$ sets,  $\{z_j\}$,  and $k$ other values on the $\{w_j\}$. This gives us two sets of $k$ values that can be ordered to be interspersed on the unit circle. Therefore, we may choose a Blaschke product $C$ of degree $k$ such that $C(0) = 0$ and $C$ identifies these two sets of $k$ values. Thus $C \circ D$ maps $0$ to $0$ and identifies the same two sets of points as $B$. As in the previous theorem, Theorem~\ref{thm:unique} implies that $B = C \circ D$.
\end{proof}

\subsection{Examples of elliptical and non-elliptical curves} In this section, we provide an example of a Blaschke product of degree $8$ with an elliptical Blaschke curve as well as a Blaschke product of degree $8$ with non-elliptical Blaschke curve.

\begin{example} We begin with an example of a Blaschke product $C$ of degree $8$ with an elliptical Blaschke curve. To connect this to the numerical range of a compression of the shift, we need $C(0) = 0$. Then $W(S_{C(z)/z})$ will have elliptical numerical range.

 Let $a \in \mathbb{D}$. Let $\varphi_a(z) = \frac{a-z}{1 - \overline{a}z}$ and consider  the automorphism $\varphi(z) = \varphi_a(z) \circ e^{i \pi/4} z \circ \varphi_a(z)$, which has the property that the composition $\varphi^{[8]}(z) = z$ for all $z \in \mathbb{C}$ and no lower-degree composition satisfies $\varphi^{[j]}(z) = z$ for all $z \in \mathbb{C}$. The Blaschke product $B(z) = \left(\varphi(z)\right)^8$ is degree $8$ and \cite[Corollary 11 ]{DGSV17} shows that for each $\lambda \in \mathbb{D}$, the line segments joining the points for which $B(z) = \lambda$ circumscribe an ellipse. By Theorem~\ref{thm:darboux2} we expect that the polygons produced by connecting vertices, as described in that theorem, will also have an ellipse (or point) as their envelope.
 
 To obtain a Blaschke product that maps the origin to itself, we consider the automorphism $\varphi_{\alpha}(z) := \frac{\alpha-z}{1 - \alpha z}$ where $\alpha = B(0)$, and note that since the Blaschke product $C(z) = \varphi_\alpha \circ B$ identifies the same points on the unit circle as $B$, the Blaschke product $C$ has an ellipse as Blaschke curve. It is easy to see that $C$ factors into a composition of degree-$2$ Blaschke products; for example, we may take $\left(\varphi_\alpha \circ z^2\right) \circ z^2 \circ \left(\varphi_a(z) \right)^2$. (It is shown in the same paper that the curves corresponding to skipping $2$ and $4$ points also product Poncelet ellipses, but this follows directly from Darboux's theorem as well.) When we take $a = .5$ we obtain the pictures below. Note that each of the ellipses in the family has two of the zeros of the Blaschke product as foci; see \cite{Mirman2} for more information on the location of foci. 
 The zeros of the Blaschke product $C$ are obtained using Mathematica: 

\begin{multline*}
 -0.158011 + 0.369131 i,
   0.0141808 + 0.629309 i, 0.241238 + 0.685693 i,\\
   0.401172 - 0.0169046 i, 0.42801 + 0.619984 i,
   0.555657 + 0.468632 i, 0.58342 + 0.236332 i  \end{multline*} plus one zero at zero.
   
\begin{figure}[ht]
 \includegraphics[width=1.5in]{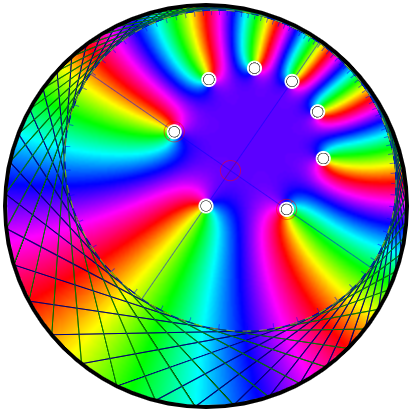}  \includegraphics[width=1.5in]{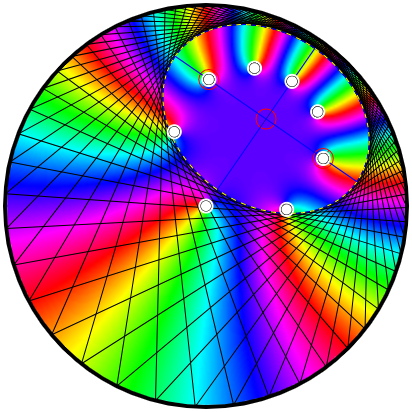}      \includegraphics[width=1.5in]{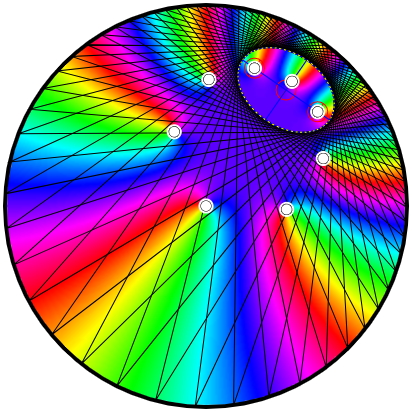}
 \caption{Degree-$8$ Blaschke product example}
\end{figure}

We say more about this in Proposition~\ref{prop:onecritical}.

\end{example}

\begin{example}\label{ex:notelliptical} We now turn to an example of a ``non-example''; that is, we give an example of a Blaschke product of degree $8$ that factors into three degree-$2$ Blaschke products but such that the Blaschke curve is {\it not} elliptical. See Figure~\ref{fig:nonexample}.

\medskip
\noindent
Let $a = .84^4$ and consider the Blaschke product 
\[B(z) = z^4 \left(\frac{z^4 - a}{1 - a z^4}\right).\] This Blaschke product has four zeros at $.84, -.84, .84 i, -.84 i$ and a zero of order $4$ at zero. It is also clear that $B = B_2 \circ B_1 \circ B_1$ where $B_1(z) = z^2$ and $B_2(z) = \frac{z (z-a)}{1 - a z}$. Due to the symmetry of the problem, we obtain vertical and horizontal tangent lines to the Blaschke curve when $B(z) = -1$. The solutions to this equation are denoted by $z_1, z_2, \ldots, z_8$. Using Mathematica we obtain the eight solutions and we are able to compute the semi-major and semi-minor axis and we find that they are both equal to $.965767$. Therefore, if this were an ellipse, the equation would be 
\[x^2 + y^2 = .965767^2.\]
By construction this circle will be tangent to the horizontal and vertical segments in the circumscribing polygon, but it must be tangent to all other sides as well. Since it must be tangent to the line segment joining $z_1$ and $z_2$, we compute the distance from the origin to this line segment. The line segment has equation:

\[ x + y = 1.22518.\] The distance from the origin to this line is 
\[\frac{1.22518}{\sqrt{2}} = .866333 \ne .965767.\] Therefore, our assumption that this is a Poncelet ellipse must be incorrect. When we connect every third point, things look quite different, see Figure~\ref{fig:not3}.

Recall \cite[Theorem~3.8]{Simanek2}, which says that if $d$ is a divisor of $n$ and $d \geq 3$, then the number of curves $C_k$, $1 \leq k \leq [n/2]$, that have the $d$-Poncelet property is $\Phi(d)/2$, where $\Phi$ is Euler’s totient function, counting the positive integers up to $d$ that are relatively prime to $d$. This implies that for $n=8$,  if $\Gamma$ is a complete Poncelet curve,  then $C_1, C_3$ are $8$-Poncelet curves, $C_2$ is a $4$-Poncelet curve and $C_4$ a 2-Poncelet curve, possibly consisting of a single point (cf. \cite[Example 3.12]{Simanek2}). Thus, Theorem~\ref{thm:Dragovic} does not apply in this setting.

\begin{figure}[ht]
 \includegraphics[width=1.5in]{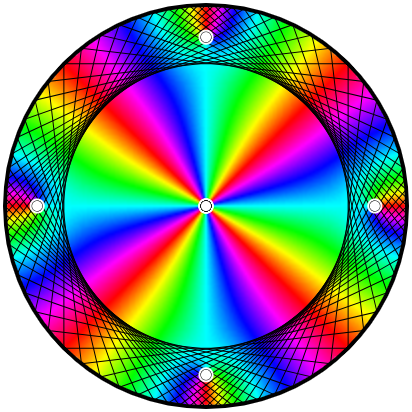}      \includegraphics[width=1.5in]{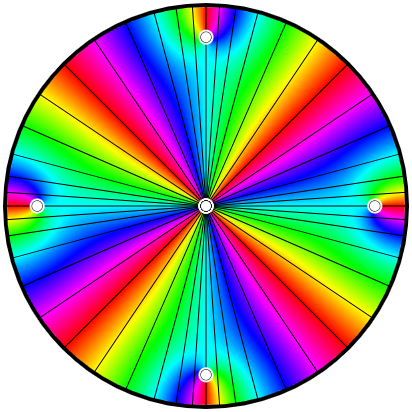}
 \caption{Degree-$8$ Blaschke product Poncelet curves (or point)}
 \label{fig:nonexample}
\end{figure}

\begin{figure}[ht]
 \includegraphics[width=1.5in]{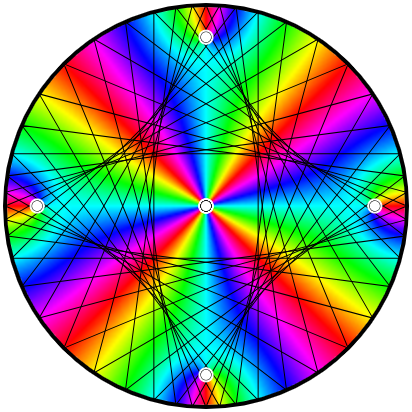}  \includegraphics[width=1.5in]{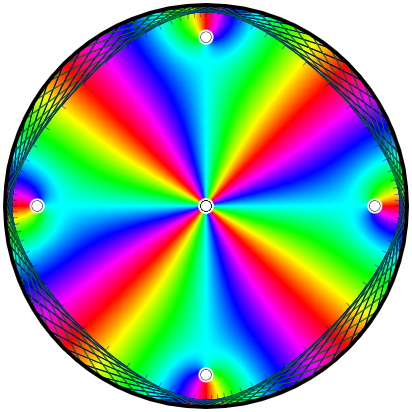} 
 \caption{Degree-$8$ Blaschke product non-conics} 
 \label{fig:not3}
\end{figure}

Now suppose that we connect every other point, which yields a convex quadrilateral.  Solving $$z^4(z^4-a)/(1-az^4)=1,$$ we see that four vertices are $1, i, -1, -i$. The other four are $\pm 1/\sqrt{2} \pm 1/\sqrt{2}$. Both of these quadrilaterals circumscribe the circle $|z| = 1/\sqrt{2}$. In this case, $D(z) = z^4$ identifies two sets of vertices of circumscribing quadrilaterals. Since two such sets of points determine such a Poncelet curve \cite{GauWu2}, $D$ identifies all points in a set of vertices and we see that the circle is the Blaschke curve associated with $D$. As a consistency check, the matrix associated with $S_{D(z)/z}$ is the $3 \times 3$ Jordan block and its numerical range is the closed disk of radius $\cos(\pi/4)$. 
\end{example}



\section{Critical Vales of Blaschke Products with Elliptical Blaschke Curve}
\label{sec:critical values}

In this section, we  study results that follow from understanding the critical values of a Blaschke product. That this is connected to composition of Blaschke products follows from results of Ritt \cite{Ritt} and Cowen \cite{Cowen}, and this will be discussed in the next section. A discussion of this also appears in the book, \cite{GarciaMashreghiRoss}. 
Note that by Theorem~\ref{thm:2}, if $B$ is of degree $2^n$ and has an elliptical Blaschke curve, then $B$ is a composition of $n$ degree-$2$ Blaschke products. 

The next theorem is useful when we count the number of critical points of a Blaschke product in $\mathbb{D}$.

\begin{theorem}[Walsh's Blaschke product theorem] \cite[p. 377]{Sheil-Small} \label{thm:Walsh} Let $B$ be a Blaschke product of degree $m$ with zeros $a_1, \ldots, a_m \in \mathbb{D}$. Then $B$ has exactly $m-1$ critical points in $\mathbb{D}$ and they lie in the convex hull of the set $\{0, a_1, \ldots, a_m\}$. The critical points of $B$ outside $\mathbb{D}$ are the 
conjugates, relative to $\mathbb{T}$, of those in $\mathbb{D}$.\end{theorem}

\subsection{Blaschke product with few critical values}
The following shows that when we have a Blaschke product of degree $2^n$ that has a Blaschke curve that is an ellipse, then we have far fewer critical values than $2^n-1$. Recall that the set of critical values is the set
\[\{w \in \mathbb{D}: w = B(z)~\mbox{and}~B^\prime(z) = 0\}.\]

\begin{example}\label{example:four} Suppose $B = C \circ D$ with $C$ and $D$ degree $2$. Then, because $B$ is degree $4$, we know that there are $3$ critical points in $\mathbb{D}$ (for example, see Theorem~\ref{thm:Walsh}) and, therefore, at most $3$ critical values. Let $z_0$ be the point with $D^\prime(z_0) = 0$ and $w_0$ the point with $C^\prime(w_0) = 0$. Let $D(z_1) = D(z_2) = w_0$. Computing $B^\prime(z) = C^\prime(D(z))D^\prime(z)$, we get critical points of $B$ at $z_0$, where $D^\prime(z_0) = 0$, and at $z_1, z_2$ where $C^\prime(D(z_1)) = C^\prime(D(z_2))  = C^\prime(w_0) = 0$. But $B(z_0)$ is one critical value and $B(z_1) = C(D(z_1)) = C(D(z_2)) = B(z_2)$ is the only possibility for the second. So $B$ has at most two critical values. We generalize this argument below in Proposition~\ref{prop:criticalvalue}. \end{example}
 




\begin{proposition}\label{prop:criticalvalue}
Let $B=B_1 \circ B_2...\circ B_n$ be a composition of $n$ Blaschke products with $\mbox{deg}~ B_j=k_j$ for $j = 1, \ldots, n$. Then $B$ has at most $\sum_{i=1}^{n} (k_i-1)$ (distinct) critical values.
\end{proposition}
\begin{proof}
We shall prove this by induction. For $n=1$, the hypothesis is true. Suppose that the hypothesis is true for $B=B_1 \circ... \circ B_{n-1}$, a composition of $n-1$ Blaschke products with $\deg B_j=k_j$ where $1\leq j \leq n-1$; that is, $B$ has a maximum of $\sum_{j=1}^{n-1} (k_j-1)$ critical values. Let  $p=\Pi_{j=1}^{n-1} k_j$ denote the degree of $B$ and consider $C=B_n\circ B$, with $\deg B_n=k_n$. Since $C'(z)=B'_n(B(z))B'(z)$, all the critical points of $B$ are critical points of $C$ that generate at most  $\sum_{j=1}^{n-1} (k_j-1)$ critical values of $C$. We know that $B_n$ has $k_n-1$ critical points and we denote them by $w_1, ..., w_{k_n-1}$. Since $B$ is a $p$ to $1$ map, there are $p$ values $z_1^{(w_j)}, ..., z_{p}^{(w_j)}$ such that $B(z_j^{(w_j)})=w_j$ for every $j$. Hence, $z_1^{(w_j)}, ..., z_{p-1}^{(w_j)}$ are critical points of $C$ for every $w_j$, but  $C(z_1^{(w_j)})=...=C(z_{k_1-1}^{(w_j)})=B_n(w_j)$, and so these critical points generate at most $k_n-1$ critical values. Then $C$ has a maximum of $(k_n-1) + \sum_{j=1}^{n-1} (k_j-1)$ critical values, as desired.
\end{proof}

\begin{corollary}
Let $B$ be a composition of $n$-degree $2$ Blaschke products. Then $B$ has at most $n$ (distinct) critical values.
\end{corollary}

The following results will be useful in this section (see also \cite{Zakeri}):

\begin{theorem}\cite{Heins}  Let $z_1, \ldots, z_d$ denote (not necessarily distinct) points in $\mathbb{D}$. Then there exists a unique finite Blaschke product of degree $d+1$ with $B(0) = 0$, $B(1) = 1$ and critical points at $z_j$. \end{theorem}

\begin{corollary}\label{cor:unique} Two proper holomorphic maps $f, g: \mathbb{D} \to \mathbb{D}$ have the same critical points, counted with multiplicity, if and only if $f = \tau \circ g$ for some conformal automorphism $\tau$ of $\mathbb{D}$.\end{corollary}

The following theorem will be useful in Section~\ref{sec:monodromy} when we discuss the monodromy group associated with Blaschke products with one critical value. We have been unable to locate a reference for this result in the literature, so we present the proof below. 

\begin{theorem}\label{onecriticalvalue} Let $B$ be a Blaschke product of degree $n$ with one critical value, $w$. Then all critical points of $B$ are equal to a point $a \in \mathbb{D}$, and there exists an automorphism $\tau$ such that 
\[B(z) = \tau \circ \left(\frac{z-a}{1-\overline{a}z}\right)^{n}.\] \end{theorem}

\begin{proof}

There exist  critical points $z_1, \ldots, z_{n-1}$ such that $B(z_j) = w$. There are $n$ points (possibly the same) satisfying $B(z) = w$, so we let these be denoted $z_0$, and $z_j$ for $j = 1, \ldots, n -1$. Let 
\[C(z) = \frac{B(z) - w}{1-\overline{w}B(z)} =  \left(\frac{z-z_0}{1-\overline{z_0}z}\right)\left(\frac{z-z_1}{1-\overline{z_1}z}\right)\cdots\left(\frac{z-z_{n-1}}{1-\overline{z_{n-1}}z}\right).\] 

Now,
\[C^\prime(z) = \frac{(1 - |w|^2)B^\prime(z)}{(1-\overline{w} B(z))^2},\] so $C^\prime(z_j) = 0$ for $j = 1, \ldots, n-1$.
So, $C$ has the same critical points as $B$ and therefore each zero, $z_1, \ldots, z_{n-1}$, is a zero of order greater than $1$, and we may
write $C$ as \[C(z) = \left(\frac{z-z_{m_1}}{1-\overline{z_{m_1}}z}\right)^{j_1}\cdots \left(\frac{z-z_{m_{l}}}{1-\overline{z_{m_{l}}}z}\right)^{j_{l}},\]
where $l<n-1$ and $z_{m_1}, \ldots, z_{m_l}$ distinct. (If no $z_j = z_0$ for $j \ge 1$, then one $j_k = 1$ and $z_{m_k} = z_0$.)

Now let us write $C(z) = \left(\frac{z-z_{m_1}}{1-\overline{z_{m_1}}z}\right)^{j_1} \times D(z)$, where $D(z_{m_1}) \ne 0$. Then $C'(z)$ has a zero at $z_{m_1}$ of at least order $j_1 - 1$. Suppose it has a zero of order strictly greater than  $j_1-1$. Then
\[C'(z) = j_1  \left(\frac{z-z_{m_1}}{1-\overline{z_{m_1}}z}\right)^{j_1-1} \times D(z) +  \left(\frac{z-z_{m_1}}{1-\overline{z_{m_1}}z}\right)^{j_1} \times D^\prime(z).\] Since $C^\prime$ is assumed divisible by  $\left(\frac{z-z_{m_1}}{1-\overline{z_{m_1}}z}\right)^{j_1}$, we would have $D$ divisible by $\left(\frac{z-z_{m_1}}{1-\overline{z_{m_1}}z}\right)$. But this is impossible because $D(z_{m_1}) \ne 0$.

Applying this argument to each factor involving $z_{m_k}$ for $k = 1, \ldots,l$, we see that each such $z_{m_k}$ can contribute at most $j_k-1$ critical points, so the total number of critical points that we get from the $z_{m_k}$ is $\sum_{k = 1}^l ( j_k - 1) = n-l.$ But there are $n - 1$ critical points, so we must have $l = 1$; in other words, 
all $z_j$ must be equal. Thus, $C$ has the same critical points as $\left(\frac{z - z_1}{1-\overline{z_1}z}\right)^n$ and by Corollary~\ref{cor:unique}, there is an automorphism $\tau_1$ such that \[C(z)  = \tau_1 \circ \left(\frac{z-z_{1}}{1-\overline{z_{1}}z}\right)^{n}.\]  Letting $\tau_w(z) = \frac{z - w}{1 - \overline{w}z}$, we have $C = \tau_w \circ B$. Thus, $B = \tau_w^{-1} \circ C$ and \[B(z) = \tau_w^{-1} \circ C(z) = \tau \circ \left(\frac{z-z_{1}}{1-\overline{z_1}z}\right)^{n},\]
with $\tau=\tau_w^{-1}\circ \tau_1.$

\end{proof}

\begin{proposition}\label{prop:onecritical}
Let $B$ be a Blaschke product of degree $n$ with one critical value. Then $B$ has an elliptic Blaschke curve.
\end{proposition}
\begin{proof}
By  Theorem 4.7, $B$ can be written as 
$$B(z)=\tau \circ \left( \frac{z-a}{1-\bar{a}z} \right)^n,$$
where $\tau$ is an automorphism. But by \cite[Corollary 10]{DGSV17} a Blaschke product of the form $\left( \frac{z-a}{1-\bar{a}z} \right)^n$ has an elliptic Blaschke curve, and this does not change by composing with automorphisms.
\end{proof}

\begin{corollary}Suppose $B$ is a Blaschke product of degree $n = p_1 p_2 \ldots p_m$ with one critical value. Then $B$ can be factored in any order as a composition of $m$ Blaschke products of degree $p_1, \ldots, p_m$. \end{corollary}

\begin{proof} This follows from the form of $B$ in Theorem~\ref{onecriticalvalue}.  \end{proof}


\section{The monodromy group and compositions of Blaschke products}\label{sec:monodromy}

We begin by considering the following from Cowen's paper \cite{Cowen}, \cite[Chapter 9]{{GarciaMashreghiRoss}}, or \cite{Wegert}. This is closely related to a theorem of Ritt.

\subsection{The decompositions of Ritt and Cowen}
The decompositions Ritt and Cowen require consideration of the critical values and a normalization of the Blaschke product. In general, if a Blaschke product has degree $n$, there are $n-1$ critical points (as is the case for polynomials) and at most $n-1$ critical values. However, as we have seen in Section~\ref{sec:critical values}, when the Blaschke product is a composition, there are fewer critical values. Following Cowen, we say that a finite Blaschke product is normalized if $B(0) = 0$, $B^\prime(0) > 0$, and $B(a) = 0$ implies that $B^\prime(a) \ne 0$.  Given a Blaschke product $B$ it is always possible to find $\alpha, \beta \in \mathbb{D}$ and $\lambda \in \mathbb{T}$ so that 
\[\lambda \varphi_\alpha \circ B \circ \varphi_\beta\] is in normalized form; see \cite[Proposition 9.2.6]{GarciaMashreghiRoss} for details.  Let $S$ denote the set of critical points, so that $B(S)$ denotes the critical values. The oriented closed loops in $\mathbb{D} \setminus B(S)$ based at the point $0$ form a group. Given two loops $\gamma$, $\delta$, the product $\gamma\cdot\delta$ is obtained by ``gluing'' them; that is, since they both start and end at zero, we begin by following $\delta$ and then continue with $\gamma$.  We consider the homotopy classes, recalling that loops are homotopy equivalent if one can be deformed to another in $\mathbb{D} \setminus B(S)$.  Since we assume that $0$ is not a critical point, $B^{-1}$ has $n$ branches at $0$ that will be denoted by $g_1, g_2, \ldots, g_n$, where $g_1(0) = 0$. We let $G_B$ be the group associated with $B$ that consists of the set of permutations of $\{g_1, \ldots, g_n\}$ induced by the loops in $\mathbb{D} \setminus B(S)$ based at $0$. The connection to composition (or, more precisely, decomposition) is described in Theorem~\ref{thm:RittCowen}. We say that a group $G$ respects a partition $\mathcal{P}$ if for each $g \in G$ and $P \in \mathcal{P}$, there exists $P^\prime \in \mathcal{P}$ so that $g P \subset P^\prime$. If $G$ respects a partition, then each element of the partition will have the same cardinality and this is called the order of $\mathcal{P}$.

For a given Blaschke product  $B$, we consider the set, $B(S)$, of critical values of $B$ and by $L_B$ we mean the set of continuous curves in $\mathbb{D} \setminus B(S)$ for which $\gamma(0)=\gamma(1)=0$.   Cowen showed \cite{Cowen}  that the monodromy group   $$G_B := \{ \gamma^*: \,\, \gamma\in L_B\}$$ can be computed from a given Blaschke product  $B$ and its local inverses; the precise statement appears in Theorem~\ref{thm:RittCowen}. He states that  if one knows all  of the normal subgroups of the  $G_B $ then  one can construct all possible non-trivial  compositional factorization  of $B$, but that this association of normal subgroups and compositions is more complicated than ``one would hope.'' Thus, we focus on the generators, rather than the group itself.

\begin{theorem}[Ritt, Cowen] \label{thm:RittCowen} Let $B$ be a finite normalized Blaschke product. If $\mathcal{P}$ is a partition of the set of branches of $B^{-1}$ at $0, \{g_1, g_2, \ldots, g_n\}$ that $G_B$ respects, then there are finite Blaschke product $J_{\mathcal{P}}$ and $b_{\mathcal{P}}$ with the order of $b_{\mathcal{P}}$ the same as the order of $\mathcal{P}$ so that $B = J_{{\mathcal{P}}} \circ b_{{\mathcal{P}}}$. Conversely, if $J$ and $b$ are finite Blaschke products so that $B = J \circ b$ then there is a partition $\mathcal{P}_b$ of the set of branches at $0$ that $G_B$ respects such that the order of $\mathcal{P}_b$ is the same as the order of $b$.
\end{theorem}

Cowen makes no claim about the equivalence of factorizations, though that is discussed in \cite{Ng13}.  Obviously, we may write $B = (C \circ \varphi_a) \circ (\varphi_a^{-1} \circ D)$ where $\varphi_a$ is an automorphism. But what is perhaps less obvious is that the degrees of the decompositions may vary. For example, $\left(\frac{z - a}{1 - \overline{a} z}\right)^6 = z^2 \circ \left(\frac{z - a}{1 - \overline{a} z}\right)^3 = z^3 \circ \left(\frac{z - a}{1 - \overline{a} z}\right)^2$. However, there is a notion of length for a Blaschke product (which requires factoring into prime factors) and the length in this case is an invariant under such factorizations, \cite[p. 24]{Ng13}.  Decompositions were also considered in \cite{Wegert}, where a method to visualize the monodromy group is presented. We give a sense of the main ideas here.  

\subsection{Visualizing the monodromy group}
Following Wegert, we consider the basins of attraction, $Z_k$, for the Blaschke product and their images, $D_k$. Recall that the basins are simply connected and their boundaries are formed by the stable manifolds $S_j$ of the critical values and arcs on the unit circle. Then $B(S_j)$ is a radial segment that has one endpoint at the critical value and $B$ maps each basin $Z_k$ in a one-to-one fashion onto slit disks $D_k = \mathbb{D} \setminus  R_k$, where $R_k$ denotes the union of all radial slits $B(S_j)$ of stable manifolds $S_j$ that belong to the boundary of $Z_k$. The Riemann surface is obtained from the $D_k$ as follows: two slit disks are connected if $Z_k$ and $Z_j$ have a common boundary component; we glue $D_k$ and $D_j$ along the image of the common component; that is, along the slit.
Because Wegert's description depends on the phase plot, it is not possible to distinguish points that are sent to values with the same argument. To handle this problem, for the following, we assume that if the Blaschke product maps zero to zero, has simple zeros, and has the property that for two critical points $z_1, z_2$, if $B(z_1)/B(z_2) > 0$, then $B(z_1) = B(z_2)$ and we will say that the Blaschke product is regularized. A Blaschke product $B$ can be regularized by composing with a disk automorphism, $\varphi \circ B$.

\begin{example}
Consider a degree-$8$ Blaschke product $B$. In visualizing what is happening using Wegert's method, due to the coloring, we assume the Blaschke product is regularized. To visualize the monodromy group, we consider loops (one from each homotopy equivalence class) that encircle the critical value exactly once. Wegert's idea is to show how the generators of the monodromy group can be read off the phase plot.

\begin{figure}[ht]
      \includegraphics[width=1.5in]{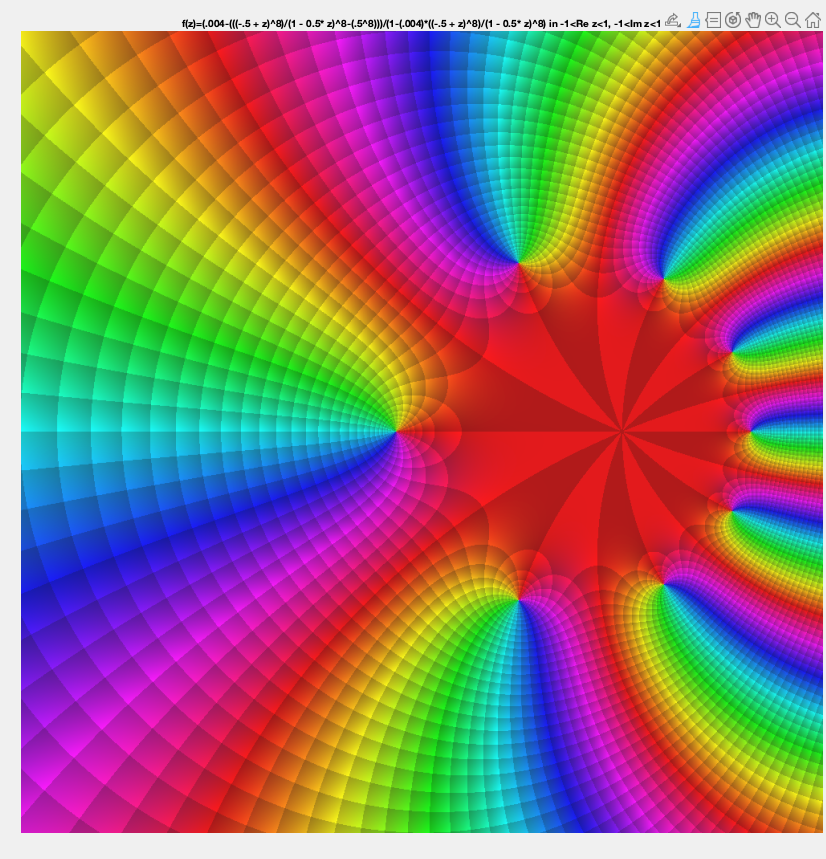}\hspace{.5in} \includegraphics[width=1.5in]{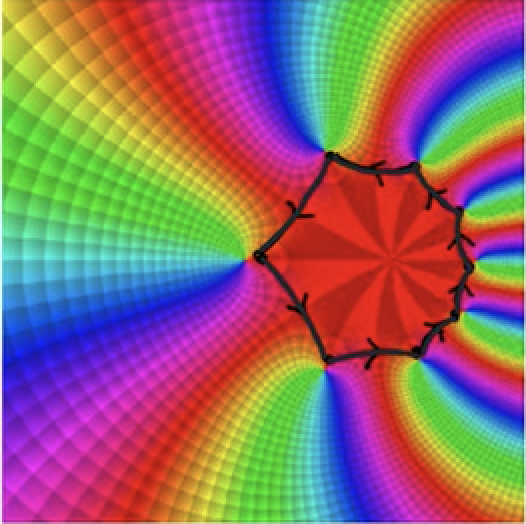}
      \caption{Blaschke product tiling and a possible generator}
            \label{degree8a}
\end{figure}

Consider a Blaschke product with a  $0$ at $0$ and seven other zeros. Figure~\ref{degree8a} is obtained via a coloring and tiling of the plane. It shows the pullback of the plane under $B$ or the phase plot of $B$.  There are seven critical points counted according to multiplicity, eight zeros, and one critical value in the red region. On a plot such as the one in this example, one can spot a critical point in the grid as a point where the grid does not form a square; that is, where the function is not conformal. These tiles, with more than four vertices, are referred to as the exceptional tiles. In this case, we see the critical value in red.  To compute the monodromy group, one needs to find the eight zeros. They are generally easy to spot because they are the places where all the colors come together. Thus, we see a critical point surrounded by $8$ zeros. Now a loop in the plane either circles the critical value or misses it. As in the particular degree-$4$ case in \cite{Cowen} (see also \cite{GarciaMashreghiRoss}), not circling the critical value corresponds to the identity map in the monodromy group. If the map circles the critical value, then when we compute the final element of the continuation, each zero will move to the next one and we obtain the generator $(12345678)$ of the cyclic group on $8$ elements. Essentially, zeros that are associated to critical points corresponding to the same critical value are moved simultaneously.

\end{example}

\subsection{Computing monodromy groups}
To illustrate this method and for later reference, we provide a detailed proof that generalizes an example of Cowen. We use his construction for Blaschke products here. We note that this also follows immediately from an observation in \cite[ p. 970]{Wegert} that generators of $G_B$ are in a one-to-one correspondence with the critical values of $B$. 

\begin{proposition}\label{thm:cyclic} Let $B$ be a normalized Blaschke product of degree $n$ with one critical value. Then the monodromy group associated with $B$ is a cyclic group of order $n$. \end{proposition}

\begin{figure}[!tbp]
  \begin{minipage}[b]{0.4\textwidth}
    \includegraphics[width=\textwidth]{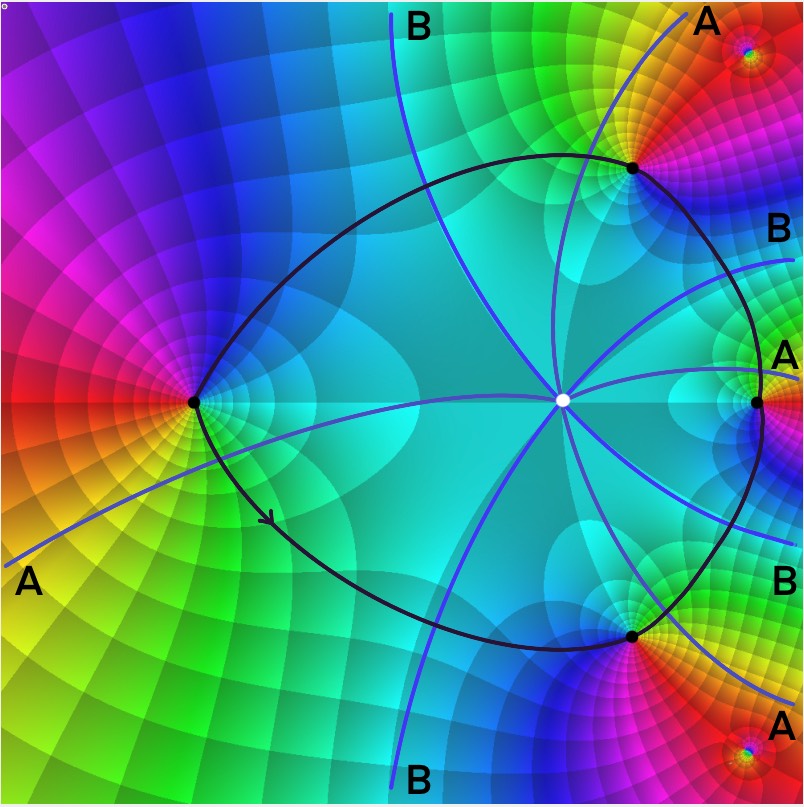}
    \caption{Inverse image.}
    \label{fig:inverseimage}
  \end{minipage}
  \hfill
  \begin{minipage}[b]{0.4\textwidth}
    \includegraphics[width=\textwidth]{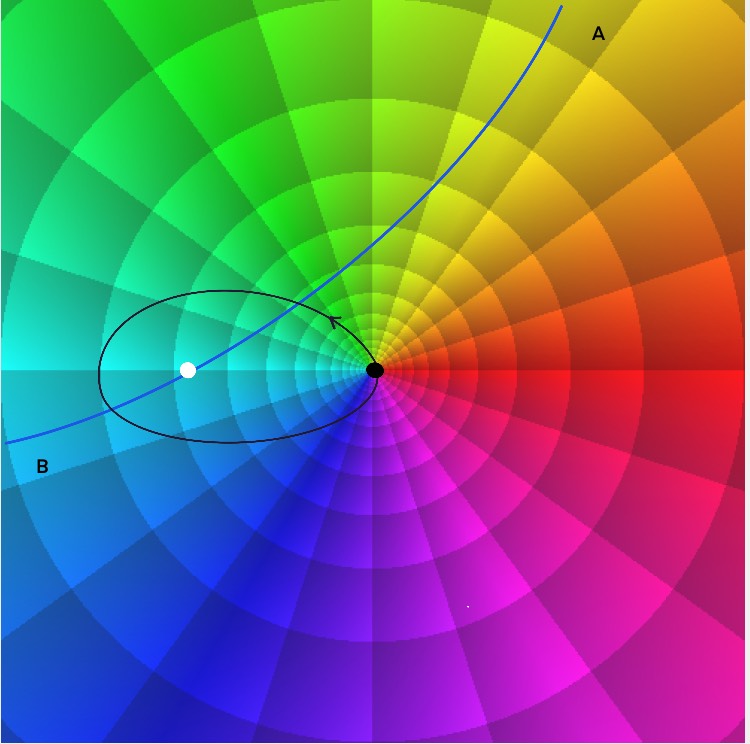}
    \caption{Loop.}
        \label{fig:onevalue}
  \end{minipage}
\end{figure}

\begin{proof}The proof is illustrated in Figure~\ref{fig:inverseimage} and Figure~\ref{fig:onevalue}. For each Blaschke product with a single critical value, we know from Theorem~\ref{onecriticalvalue} that there will be one critical point of order $n-1$ in $\mathbb{D}$. Draw a path starting at $A$, through the critical value and ending at $B$ (picture on the right). If the Blaschke product is of degree $n$, the inverse image will have $n$ curves, each passing through the critical point (picture on the left). Note that  because the argument of the Blaschke product is increasing, the inverse images of the points $A$ and $B$ are interlaced. The inverse image of a loop based at the origin will begin at a zero of the Blaschke product and, if oriented counterclockwise, will always pass through the curve associated with $A$ (the purple curve) between the critical point and the inverse image of $A$. The inverse image  must  then pass between the critical point and $B$ and it must end at a zero. Thus, it has moved from a zero $z_1$ to a zero $z_2$. This will be repeated until the curve returns to $z_1$. Thus, the permutation associated with this is $(123\ldots n)$. If the loop does not contain the critical value, then we obtain the identity. So the monodromy group is the cyclic group on $n$ elements.  
\end{proof}

We begin with more discussion of Blaschke products of degree $4$. These will also serve as the first step of the induction in Theorem~\ref{thm:mono}

\bigskip
\begin{example}\label{ex:degree4} Let $B$ be a normalized Blaschke product of degree $4$ with $B = C \circ D$ with $B(0) = 0$. As in \cite[Proposition 2.1]{DGSSV15} we may assume that $B(0) = C(0) = D(0) = 0$, $C$ and $D$ degree $2$, and $B, C,$ and $D$  regularized. In \cite[Section 9.4]{GarciaMashreghiRoss} two (modified) examples of Cowen's Blaschke products of degree $4$ that are compositions are presented. In one case, the group is computed to be the cyclic group of order $4$ generated by $(1234)$. This example has one critical value. The second example has two critical values. Here, the monodromy group is shown to be the dihedral group on $4$ elements (order $8$). We show that this works in general, for regularized Blaschke products of degree $4$ that are compositions of two (regularized) degree-$2$ Blaschke products and have two distinct critical values. We consider two cases.

\begin{enumerate}
\item \label{ex:4.1} If $B$ has one critical value, then by Proposition~\ref{thm:cyclic}, the monodromy group is cyclic and generated by the permutation $(1234)$.
\item \label{ex:4.2} Suppose that $B$ has two critical values. Note that since we assume $D(0) = 0$ we must have $D(z) = z (a - z)/(1 - \overline{a} z) = z \varphi_a(z)$ for some $a \in \mathbb{D}$.  Since we assume that $B$ has simple zeros, we may assume that $a \ne 0$. So, $D(\varphi_a(z)) = \varphi_a(z)\left(\varphi_a(\varphi_a(z))\right) = z \varphi_a(z) = D(z)$. Therefore, 
\[B(\varphi_a(z)) = C \circ D(\varphi_a(z)) = C \circ D(z) = B(z).\] So, \[B^\prime(\varphi_a(z)) \varphi_a^\prime(z) = B^\prime(z).\] But $\varphi_a$ has no critical point in $\mathbb{D}$ so we have  $B^\prime(z) = 0$ if and only if $B^\prime(\varphi_a(z)) = 0$.  

Since the set of critical points $\{z_1, z_2, z_3\}$ must be invariant under $\varphi_a$ and $\varphi_a$ is self-inversive, either $\varphi_a(z_3) = z_j$ for $j = 1,2$ or $\varphi_a(z_3) = z_3$. If the critical points are distinct, then we may assume that $z_2 = \varphi_a(z_1)$.  Because we assume the points are distinct, we can only have $\varphi_a(z_3) = z_3$. A computation shows that, since $a \ne 0$, we have $z_3 = \frac{1-\sqrt{1-|a|^2}}{\overline{a}} = a_\star$. 

If the critical points are not distinct, then two must be equal, say $z_1 = z_2$. If $\varphi_a(z_1) \ne z_1$, then the third must be $\varphi_a(z_1) = z_3$. But then $B(z_1) = B(z_2) = B(\varphi_a(z_1)) = B(z_3)$ and there is only one critical value, which is the case that we have already handled. Therefore $z_1 = z_2 = a_\star$ and either $\varphi_a(z_3) = z_3$ or $\varphi_a(z_3) = z_1$. Either way, all three points must be equal and $B$ would have only one critical value.

We now turn to the monodromy group in this case. We know that $B$ has three critical points $z_1, \varphi_a(z_1),$ and $a_\star$ and four zeros. In this case, a loop circling a critical value corresponding to one critical point (the critical point of $D$) will yield a transposition. We assume, upon re-numbering, that it is $(13)$. We know that there are just two critical values and therefore two generators. Each zero of $D$ will remain when we compose with $C$ and we add two more zeros, one to each basin, and two more critical points. Thus, a loop circling one critical value corresponding to two critical points will move points simultaneously passing about the other critical points and therefore will yield a generator that is a product of transpositions, say $(12)(34)$, while a generator circling both will yield $(1234)$, which is the product of these two. We may choose either one as our second generator. Thus we can replace these three generators with just two, namely
$(13)$ and $(12)(34)$, which is the number of critical values.
\end{enumerate}
This second set of generators yields a group that is isomorphic to the wreath product of two cyclic groups of order two; that is, the group, $\mathbb{Z}_2 \wr \mathbb{Z}_2$. We say more about this below. \end{example}

\subsection{Wreath products and trees.}  The theory of wreath products is often explained  by thinking of them as groups acting on a finite rooted tree. We are grateful to Peter Brooksbank for providing this background on wreath products.

 Let $\mathcal{T}= \mathcal{T}_k$ denote a binary tree of height $k$, where $k \geq 1$ and let $n=2^k$.  Setting $\Omega=\Omega_k = \{1,\dots, 2^k \}$, one can label nodes and leaves of the tree as follows: start labelling the root node by $\Omega$ and rest of the nodes of $\mathcal{T}$ with subsets of $\Omega$. If a node at level $0\leq j <k$ is labelled with $\{ m+1,m+2,\dots,  m+2^{k-j} \}$ for some integer $m$, one can label its left child at the level $j+1$ with $\{m+1, \dots, m+2^{k-(j+1)}\}$ and the right child with $\{ m+2^{k-(j+1)} +1 , \dots, m+2^{k-j}\}$.  In this manner the leaves of the tree are labelled from left to right  with $\{1\},\dots,\{n\}$ or, equivalently, $1,\dots,n $. 
 Let $\Gamma=\Gamma_k$ and Sym$(\Omega)$ denote the group of automorphisms of $\mathcal{T}$ and symmetric group on $\Omega$, respectively. One can identify $\Gamma$ with its image in Sym$(\Omega)$ by noticing that  the action of $\Gamma$ on $\Omega$ gives a faithful representation $\Gamma \to $ Sym$(\Omega)$.  
(See pages 45-50 in \cite{DM} for the definition of wreath product). Additionally, since $\Gamma$ permutes the nodes at each level, the block labels are permuted in the action at each level making it possible to view $\Gamma$ and its iterated wreath product structure.
 
 Consider the following elements of Sym$(\Omega)$ :
 
 $$\sigma_k:= (12), \,\sigma_{k - 1} = (13)(24),  \dots ,\sigma_{1} :=\prod_{l=1}^{2^{k-1} }(l (l+2^{k-1})).$$
 Then each $\sigma_j$ is an automorphism of the labelled tree $\mathcal{T}$; in fact $\sigma_j$ is an automorphism of the leftmost subtree rooted on level $j-1$.
 
\begin{figure}[ht]
\includegraphics[width=3in]{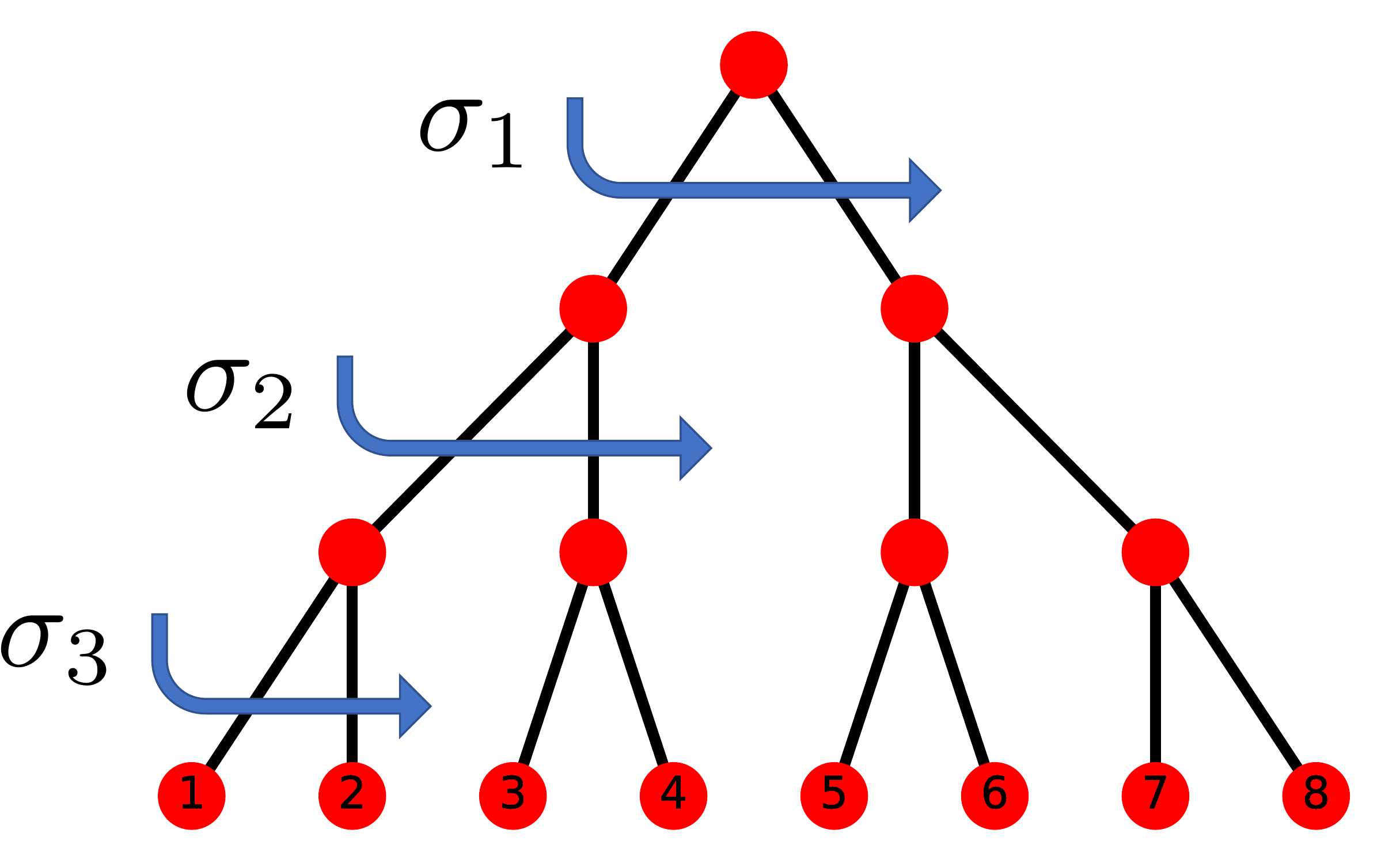}
\caption{Tree for $k = 3$}
\label{fig:tree}
\end{figure}

For example, the tree for $k = 3$ is presented in Figure~\ref{fig:tree}. The group is a semi-direct product generated by $\sigma_3 = (12), \sigma_2 = (13)(24), \sigma_1 = (15)(26)(37)(48)$ and has $128 = 2^{1 + 2 + 4}$ elements. 

 Next we present a well-known result  on three basic properties of the group of automorphisms of $\mathcal{T}$ (see, for example, \cite[p. 140]{Rotman} for \eqref{item:2} and \eqref{item:3}). In order to make the paper more accessible, we  include a short proof.
  
 \begin{proposition}
 The following hold true for the group of automorphisms of $\mathcal{T}$:
 
 \begin{enumerate}
 \item \label{item:1} $\Gamma =\langle \sigma_1, \dots \sigma_k\rangle$;
 \item \label{item:2} $\Gamma$  is an iterated wreath product of $k$ cyclic groups of order $2$;
 \item \label{item:3} $\Gamma$ is a Sylow $2$-subgroup of Sym$(\Omega)$.
 \end{enumerate}
 \end{proposition}
 
 \begin{proof}
 
The proof of (\ref{item:1})  is done by induction on $k$, which is the height of the tree. Since the result is clear for  $k=1$, let us assume $k>1$ and assume  the result holds true for trees of height less than $k$. Notice that each automorphism 
 $\alpha \in \Gamma$ permutes the two nodes on level $1$; these two nodes are the children of the root node.  Observe that $\alpha$ either fixes both children or interchanges them and we want to prove that $\alpha\in \langle \sigma_1, \dots \sigma_k \rangle$.
 If $\alpha$ interchanges the children, replace $\alpha$ with $\alpha\sigma_1$ so that $\alpha$ now fixes both children. We may now write  \[ \alpha= \beta \gamma=\gamma \beta,\] where $\beta$ is an automorphism of  the left subtree that is the identity on $\{ \frac{n}{2}+1,\dots, n\}$  and $\gamma$ is an automorphism  of the right subtree that is the identity on $\{1,\dots, \frac{n}{2}\}$.

By induction, $\beta\in \langle\sigma_2, \dots, \sigma_k\rangle$.  Additionally,  $\sigma_1\gamma\sigma_1$ is an automorphism of the left subtree and it is the identity on the right subtree. Therefore $\sigma_1 \gamma \sigma_1 \in \langle \sigma_2, \dots, \sigma_k \rangle$. Hence,
 $$ \alpha= \beta \gamma \in   \langle \sigma_1, \dots \sigma_k \rangle.$$
 That $\Gamma$  is  an iterated wreath product of $k$ cyclic groups of order $2$ follows from the construction. To check that   $\Gamma$ is a Sylow $2$-subgroup of Sym$(\Omega)$, we have only to compute the order of $\Gamma$, which is
 $$ |\Gamma| = 2^{1+2+\dots+2^{k-1}}= 2^{2^{k}-1}.$$
\end{proof}

\begin{remark} Consider a Blaschke product of the form $B(z)= z^2 \left(\frac{a-z}{1-a z}\right)^2$ with $a \in [0,1] \setminus \{0\}$. After normalizing to obtain a Blaschke product with two critical values, it follows from Example~\ref{ex:degree4} that the monodromy group is the dihedral group $D_8$.   Thus, its order is  $8$ and  this group is non-abelian.  It has two subgroups,  $H=\{1,s\}$ of order 2 and $K=\{1,r,r^2,r^3\} $ of order 4.  Since all groups of order at most $5$ are abelian, we know both $H$ and $K$ are abelian. Since direct product of two abelian groups is abelian, $D_8$ cannot be a direct product of its subgroups (this shows that \cite[Theorem 4]{Wegert} is not correct as stated). However,  $D_8$  is a semidirect product of  the same subgroups $H$ (reflection over the diagonal) and $K$ (rotation by an angle $ \pi/2$).  Furthermore, since every subgroup of index $2$  is a normal subgroup, $K$ is a normal subgroup.  Since the wreath product is a special combination of two groups based on the semidirect product, it is also possible to express  $D_8$ as the wreath product of $\mathbb{Z} /2\mathbb{Z}\, \wr\,  \mathbb{Z} /2\mathbb{Z}$.
By the definition of the wreath product this means:
$$(\mathbb{Z} /2\mathbb{Z}  \times \mathbb{Z} /2\mathbb{Z})  \rtimes  \mathbb{Z} /2\mathbb{Z}.$$
 If we think of $D_8$  as the group of automorphism of the square, the term $(\mathbb{Z} /2\mathbb{Z}  \times \mathbb{Z} /2\mathbb{Z})  $  corresponds to the swapping the sides of the square  (rotations) and the term at the end $ \mathbb{Z} /2\mathbb{Z} $ corresponds to swapping through the diagonals (reflections).
\end{remark}

Before we prove the next theorem, we note that an example from \cite{GorkinWagner} shows that  even if a Blaschke product is decomposable,  the order of the composition matters and can change the shape of the boundary of $W(S_B)$. In fact, in \cite[Theorem 10]{Simanek3}, the authors show that for a degree-$6$ Blaschke product, the numerical range being elliptical is equivalent to the corresponding Blaschke product $\widehat{B}$ having the property that it factors as  $C_1 \circ D_1$ with $C_1$ degree-$2$ and $D_1$ degree $3$ and it also factors as $D_2 \circ C_2$ with $C_2$ degree $2$ and $D_2$ degree $3$. For example, if we let $\widehat{B}_1$ be a Blaschke product of degree-$6$ such that  $\widehat{B}_1 = C_1 \circ D_1$ where 
$$C_1(z) = z  \left( \frac {z-a}{1-\overline{a}z} \right)^2\quad\mbox{and} \quad D_1(z)= z^2.$$  then one can show $W(S_{B_1})$ is an elliptic disk, where $\widehat{B}_1(z) = z B_1(z)$. Note that 
\[\widehat{B}_1(z) = z^2 \left( \frac {z^2-a}{1-\overline{a}z^2} \right)^2 = z^2 \circ z\left( \frac{z^2-a}{1-\overline{a} z^2}\right) = z\left( \frac{z-a}{1-\overline{a} z}\right)^2 \circ z^2 .\]  On the other hand, if we consider  $\widehat{B}_2= C_2 \circ D_2$ where
$$C_2(z) = z \left( \frac{z-.5}{1-.5z}\right) \quad \mbox{and}\quad D_2(z)= z^3, $$  then letting $\widehat{B}_2(z) = z B_2(z)$,  it turns out that  $W(S_{B_2})$  is not an elliptical disk.

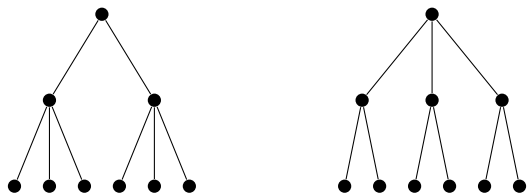
\begin{figure}
    \centering
\begin{forest}
  for tree={%
    label/.option=content,
    grow=south,
    content=,
    circle,
    fill,
    minimum size=5pt,
    inner sep=0pt,
  }
[
[[] [] [] ]
[[] [] [] ]
]
\end{forest}
\qquad \qquad
\begin{forest}
  for tree={%
    label/.option=content,
    grow=south,
    content=,
    circle,
    fill,
    minimum size=5pt,
    inner sep=0pt,
  }
[
[ [] [] ]
[ [] [] ]
[ [] [] ]
]
\end{forest}
    \caption{The tree for $S_2 \wr S_3$ (left) and $S_3 \wr S_2$ (right).}
        \label{fig:my_label}
\end{figure}

As discussed above, the groups that we consider will have monodromy groups that can be represented by trees. Here we see the same effect on order: Suppose we have a Blaschke product of degree $6$ decomposed as  $C$ of degree $2$ and  $D$  of degree $3$.   The tree representation of the wreath product  $S_2 \wr S_3$  of the monodromy group  $G_B$  with $B=C\circ D$ is quite different from the tree representation of the  wreath product of $S_3\wr S_2$ of the group $G_B$ when $B=D\circ C$ as shown in the Figure~\ref{fig:my_label}.


\begin{theorem}\label{thm:mono} Let $B$ be a regularized Blaschke product with distinct critical points and $n$ distinct critical values that is the composition of $n$ degree-$2$ regularized Blaschke products; that is, $B = B_j \circ \cdots \circ B_1$ with $B_j$ regularized for each $j$. Further assume that $B_j(0) = 0$ for $j = 1, \ldots, n$. Then the monodromy group associated with $B$ is the wreath product of $n$ cyclic groups of order $2$, or $\underbrace{\mathbb{Z}_2 \wr \mathbb{Z}_2 \cdots \wr \mathbb{Z}_2}_{n~ times}$.
\end{theorem}

The general form of an element in such a wreath product, depending on the numbering, would be
\[\prod_1^1 (i ( i+1)) = (12), \, \prod_{i = 1}^2 (i (i + 2)) = (13)(24),\,  \ldots, \prod_{i = 1}^{2^{n-1}}(i (i+2^{n-1})).\]

\begin{proof}
Note that the assumption that $B_j(0) = 0$ for each $j$ implies that $B_j \circ \cdots \circ B_1$ is a subfactor of $B_{j + 1} \circ \cdots \circ B_1$ for each $j$. Thus, the zeros of $B_j \circ \ldots \circ B_1$ are a subset of the zeros of $B_{j+1} \circ B_j \circ \dots \circ B_1$ and all such zeros are simple. We have shown in Example~\ref{ex:degree4} that the result is true for the composition of two degree-$2$ Blaschke products. So suppose that it is true whenever we have the composition of $n-1$ degree-$2$ Blaschke products. Consider $B_n \circ B_{n-1} \circ \ldots \circ B_1$. Now by our assumption, each $B_j$ must give rise to a distinct critical value. Notice that $B_1$ will correspond to one critical point, $B_2$ to two additional critical points, and, in general, $B_j$ will add $2^{j-1}$ critical points to those obtained from $B_{j-1} \circ \cdots \circ B_1$. Considering the loops that do not circle the critical value associated with $B_n$ we obtain, by induction, $n-1$ generators that yield the group $\underbrace{\mathbb{Z}_2 \wr \mathbb{Z}_2 \cdots \wr \mathbb{Z}_2}_{n - 1~ \mbox{times}}$. If a loop circles the new critical value that we obtain from $B_n$, then with an appropriate numbering, all critical points that are added are associated with the same critical value and so the action on each will be simultaneous. Note that the additional composition with $B_n$ adds a zero to each basin and a critical point, which is the one that will be encircled. Thus each new zero is paired with a zero that was in the previous basin and we obtain a generator $(1(2^{n-1}+1))\cdots(2^{n-1}2^n)$ as our generator. Since the number of generators is the same as the number of critical values, we now have the complete set of generators. Thinking of our wreath product as a group acting on a finite rooted tree, \cite[Theorem 2.1.6]{CST} completes the proof of the theorem. \end{proof}

Remark: As in Example~\ref{ex:4.2}, this is not the only generator one can choose. With proper ordering, there is a cycle of length $n$ that serves as generator.\\

\section{Groups of invariants}\label{sec:group of invariants}

Given a finite Blaschke product $B$ of degree $n$, denote the set of continuous functions $u : \mathbb{T} \to \mathbb{T}$ by $C(\mathbb{T})$ and consider the group of invariants of $u$ defined by \[\mathcal{G}_B = \{u \in C(\mathbb{T}): B \circ u = B\}.\]  This set $\mathcal{G}_B$ is a semigroup under the operation of composition; that is, the composition of two functions in $C(\mathbb{T})$ is again in the set, the identity is in the set and the operation is associative. We also note that in order to be in the group of invariants a continuous function $u$ must be a bijective mapping of the circle (see \cite[Lemma 4.1]{CGP12}, where these results are extended to infinite Blaschke products) and, for later reference, we note that the argument of $B$ (appropriately chosen) is increasing on the unit circle and therefore the zeros of $B(z) - \lambda$ are simple for every point $\lambda \in \mathbb{T}$. (This is well known; for a reference see for example, \cite[Remark 2.1]{CGP12}).  Thus, each element in $\mathcal{G}_B$ has an inverse. Cassier and Chalendar showed \cite{CassierChalendar} that 
$\mathcal{G}_B$ is a cyclic group of order $n$. In this section, we consider composition and its relation to the group of invariants. One such theorem is given below.

\begin{theorem} \cite[Theorem 5.13]{DGSSV15} A Blaschke product $B$ of degree $n = m k$ with $m > 1$ is a composition of two nontrivial Blaschke products if and only if there exists a Blaschke product $D$ of degree $k > 1$ such that $\mathcal{G}_D$ is generated by $g^m$ for some generator $g$ of $\mathcal{G}_B$. \end{theorem}

If $D$ has degree $2$, then $\mathcal{G}_D$ has order $2$. A generator of the group can be found using an observation of Frantz \cite{Frantz}: Let $a \in \mathbb{D}$. If $z \in \mathbb{T}$ and we consider the line through $z$ and $a$. Since this is not tangent to $\mathbb{T}$, Frantz shows that $\varphi_a(z) = (a-z)/(1 - \overline{a} z)$ is the other point of intersection of the line joining $z$ and $a$ with $\mathbb{T}$. We use this in Theorem~\ref{thm:invariants}.



For compositions of degree-$2$ Blaschke products, it's possible to say more. For example, if $C_j$ are degree-$2$ Blaschke products for $j = 1,2,3$ and $B = C_3 \circ C_2 \circ C_1$, there is a generator $w$ such that the corresponding group is $\langle w \rangle$ where $w^2 = v$ is the generator of $\mathcal{G}_{C_2 \circ C_1}$, and $w^4 = v^2 = u$ is the generator of $\mathcal{G}_{C_1}$. We extend this observation below.  


Note that we may write a factorization of $B$ as $B = C_n \circ C_{n-1} \circ \ldots \circ C_1 =  C_n \circ C_{n-1} \circ \ldots  \circ \varphi_{C_1} \circ \varphi_{C_1} \circ C_1$, and $\varphi_{C_1} \circ C_1$ has the same group of invariants as $C_1$. Thus, we may assume that $C_1(0) = 0$.

\begin{theorem}\label{thm:invariants} Let $B$ be a Blaschke product of order $2^n$ and let $g$ denote a generator of $\mathcal{G}_B$. If $B = C_n \circ C_{n-1} \circ \cdots \circ C_1$ is a composition of $n$ Blaschke products of order $2$, then the group of invariants of $C_j \circ \cdots \circ C_1$ is a normal subgroup of index $2$ of the group of invariants of $C_{j+1} \circ C_j \circ \cdots \circ C_1$ for each $j$.  If  $C_1(z) = z \left(\frac{a-z}{1-\overline{a}z}\right)$, then $g^{2^{n-1}} = \varphi_a$.
 \end{theorem}

\begin{proof} 

The first statement follows by using the fact that the group of invariants of $C_j \circ \cdots \circ C_1$ is a cyclic group of order $2^j$: Every subgroup of a cyclic group is normal and using Lagrange's theorem we conclude that the index is $2$. Therefore, the group of invariants of $C_j \circ \cdots \circ C_1$ is a normal subgroup of index two of the group of invariants of $C_{j+1} \circ C_j \circ \cdots \circ C_1$. We also know that since $\mathcal{G}_B$ is a cyclic group, if $g$ is a generator of $\mathcal{G}_B$ and $H$ is a nontrivial subgroup, then $g^m$ is a generator of $H$, where $m$ is the smallest positive integer with $g^m$ in $H$. In our case, we take the subgroup $H_j$ corresponding to the group of invariants of $C_j \circ \cdots \circ C_1$, which has order $2^j$. Thus, if $g^m$ generates $\mathcal{G}_{C_j \circ \cdots \circ C_1}$, then $m$ is the smallest integer such that $(g^m)^{2^j} = e = g^{2^n}$. So, $m = 2^{n - j}$.

To prove the remaining assertion, first suppose that $B$ has degree $2$; that is, that $n = 1$. As above, we may assume that $B(0) = 0$ and we may write $B(z) = \lambda z \left(\frac{a-z}{1 - \overline{a} z}\right)$.  We will suppose that $\lambda = 1$.

We have for every $a \in \mathbb{D}$,  $\varphi_a \ne id$ and $B \circ \varphi_a(z) = \varphi_a(z) (\varphi_a \circ \varphi_a)(z) = B(z)$. Thus $\varphi_a$ is a generator of $\mathcal{G}_B$. Note that $\varphi_a$ has no fixed point on the unit circle and maps a point of $w_1 \in \mathbb{T}$ to a second (distinct) point $w_2$ on $\mathbb{T}$. Since $B \circ \varphi_a = B$, we see that $B(w_1) = B(w_2)$. 

Now suppose that $u \ne e$ is another generator of $\mathcal{G}_{B}$, and consider two points $z_1, z_2$ that $B$ identifies. Either $u(z_1) = z_2 = \varphi_a(z_1)$ and $u(z_2) = z_1 = \varphi_a(z_2)$, or there exist $z_1, z_2$ with $u(z_1) = z_1$ and $u(z_2) = z_2$. If the latter case occurs for some $z_1$ and $z_2$, then this divides the unit circle into two arcs, $\ell_1$ and $\ell_2$ with endpoints $z_1$ and $z_2$, see Figure~\ref{arcs}. Now either $u(\ell_1) \subseteq \ell_1$ or $u(\ell_1) \subseteq \ell_2$. In the first case, all points in $\ell_1$ would be mapped to themselves under $u$ and the same is true for $u(\ell_2)$. So $u = e$. Thus, $u$ must interchange the two arcs. Choose a sequence $(w_n)$ in $\ell_2$ tending to $z_1$. Then, by the discussion above, (see also \cite{Frantz}) $u(w_n)$ is equal to the endpoint of the line segment on $\mathbb{T}$ that passes through the point $a$. Therefore $u(w_n) \to z_2$. So $u(z_1) = z_2 = \varphi_a(z_1)$. Since this is true for all points on $\mathbb{T}$, we see that $\varphi_a$ is the only generator of $\mathcal{G}_B = \mathcal{G}_{C_1}$.

\begin{figure}[ht]
 \includegraphics[width=1.5in]{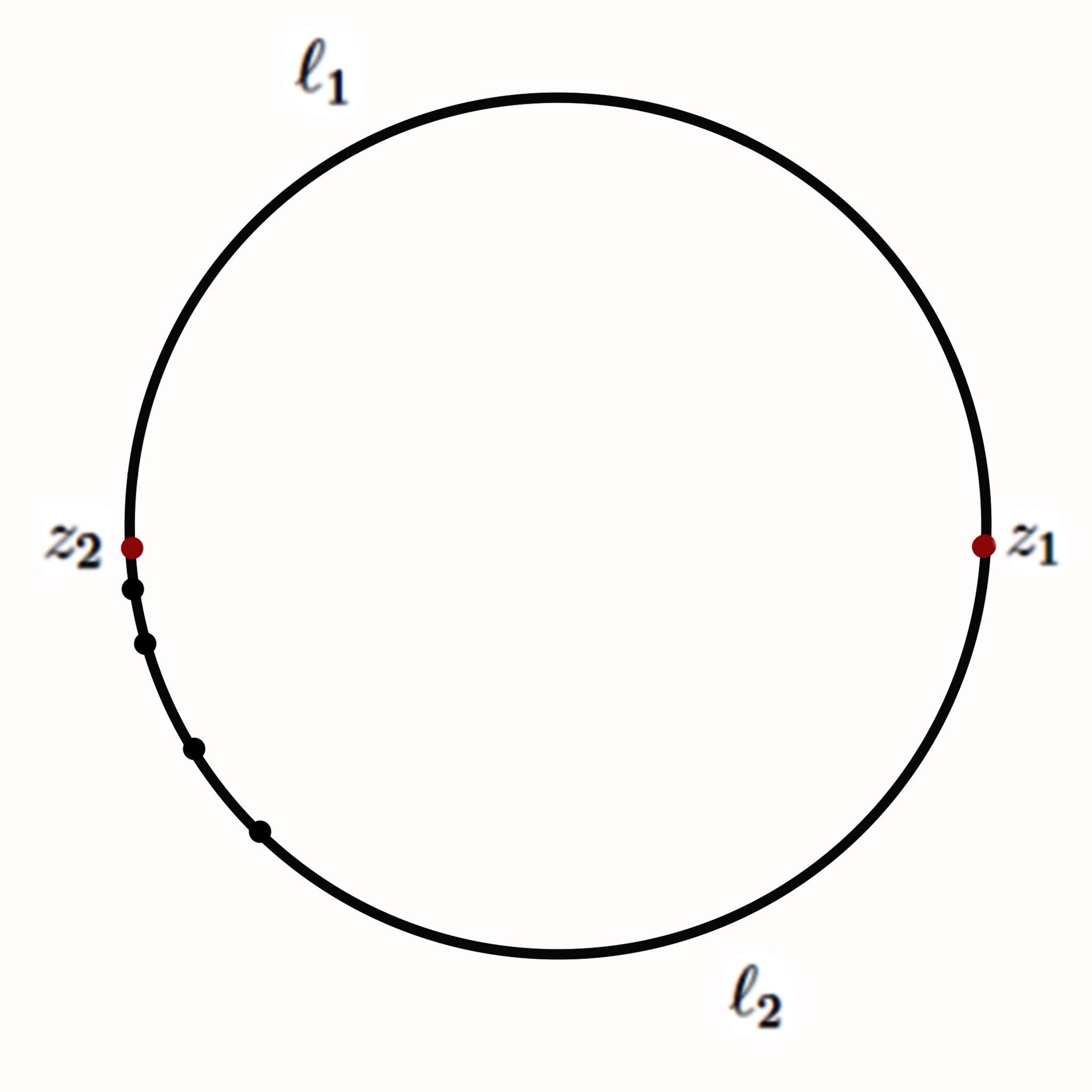}      \includegraphics[width=1.5in]{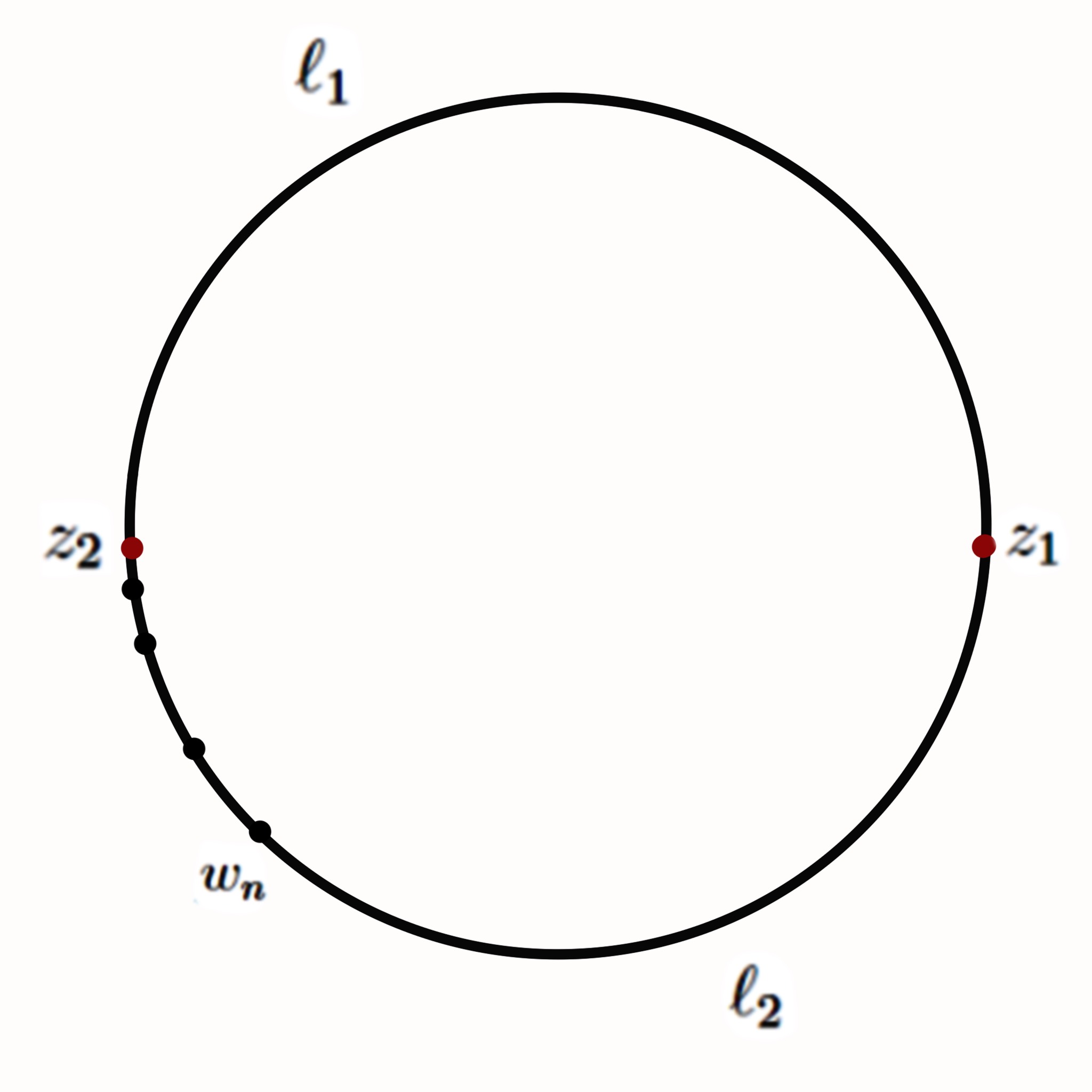}
 \caption{Arcs  $\ell_1, \ell_2$ and points $z_1, z_2$}
 \label{arcs}
\end{figure}

Now suppose that $n > 1$ and $B = C_n \circ \cdots \circ C_1$. Then by our work thus far, $\mathcal{G}_B$ has a generator $g$ of order $n$ and $g^{2^{n-1}}$ is a generator of $G_{C_1}$. But  every generator of $\mathcal{G}_{C_1}$ must identify the same points as $\varphi_a$. Thus, $g^{2^{n-1}} = \varphi_a$.

\end{proof}

\section{Final Remarks}
Mirman considered conditions for a Poncelet curve to close, the connection of Poncelet's theorem to boundaries of the numerical ranges, ways to locate the foci of the algebraic curves obtained, and a related measure in \cite{Mirman, Mirman2, Mirman98}.  A detailed look at Mirman's work on ellipses and Poncelet packages, with related examples and counterexamples, as well as results on orthogonal polynomials on the unit circle can be found in {\cite{Simanek2}. Work on composition for Blaschke products of degree $2n$ and $3n$ can be found in \cite[Corollary 5.8 and Corollary 5.10]{DGSSV15}. A description of ellipses  when the Blaschke product $\widehat{B}$  defining the Blaschke curve (that is, the boundary curve for $W(S_B)$) has degree-$6$ can be found in \cite{Simanek3}. 

Considering a Blaschke product $B$ with zeros $a_1, \ldots, a_n$, the polynomials $\Phi_n(z) = \prod_{j = 1}^n(z - a_j)$ and the reversed polynomials $\Phi_n^\star(z) = z^n \overline{\Phi_n(1/\overline{z}})$, then \[B(z) = \frac{\Phi_n(z)}{\Phi_n^\star(z)}\] and many of the results that we use and prove here can be restated in terms of the polynomials $\Phi_n$ and CMV matrices.  See \cite{Simon} for details. For example, compositions of Blaschke products can be rephrased as compositions of orthogonal polynomials, \cite[Theorem 1]{Simanek3} and, using this, many results stated here are easily restated in this setting. Connections between the theory of orthogonal polynomials on the unit circle with Poncelet's theorem, Blaschke products, and numerical ranges can be found in \cite{Simonetal}.

\section*{Acknowledgement}

The authors are grateful to Peter Brooksbank and Elias Wegert for many valuable conversations. We also thank the organizers of the ``Working Groups for Women in Operator Theory 2021'' held virtually at the Lorentz Center.

\end{document}